\documentclass[review]{elsarticle}

\usepackage{lineno,hyperref}
\modulolinenumbers[5]

\journal{~}

\newtheorem{algorithm}{Algorithm}[section]
\newtheorem{remark}{Remark}[section]
\usepackage{exscale}
\usepackage{graphicx,amsmath,bm}
\usepackage{amssymb,amscd,verbatim}
\usepackage{color}
\usepackage{float}
\usepackage{algpseudocode,algorithm,algorithmicx}

\algrenewcommand\algorithmicrequire{\textbf{Precondition:}}
\algrenewcommand\algorithmicensure{\textbf{Postcondition:}}

\usepackage{amsmath}
\usepackage{amsfonts}
\usepackage{amssymb}
\usepackage{amstext}
\usepackage{amsbsy}
\usepackage{arydshln}
\usepackage{amsmath}
\usepackage{xifthen}
\usepackage{graphicx}
\usepackage{amsmath}
\usepackage{amsfonts}
\usepackage{amssymb}
\usepackage{amstext}
\usepackage{amsbsy}
\usepackage{epsfig}
\usepackage{epsf}
\usepackage{multirow}
\usepackage{float}
\usepackage{stmaryrd} 

\usepackage{xspace}
\usepackage{url}
\usepackage{verbatim}
\usepackage{ifthen}
\usepackage{rotating}
\usepackage{mathpazo}
\usepackage{bbding}
\usepackage{booktabs}
\usepackage[bbgreekl]{mathbbol}
\usepackage{amsfonts}
\DeclareSymbolFontAlphabet{\mathbbl}{bbold}
\usepackage[noprefix]{nomencl} 
\usepackage{bbm, dsfont}

\newcommand{\sphere}{\ensuremath{\mathcal{S}}\xspace}
\newcommand*{\domain}{\ensuremath{\mathcal{D}}\xspace}
\newcommand{\unit}[1]{\,\mathrm{#1}}
\newcommand{\afluxunit}{\,\ensuremath{\unit{cm}^{-2}\unit{s}^{-1}\unit{st}^{-1}}\xspace}
\newcommand{\direction}{\ensuremath{\vec{\Omega}}\xspace}
\newcommand{\grad}{\vec{\nabla}}

\newcommand{\Sigt}[1][]{\ensuremath{\Sigma_{\xslabel[#1]{t}}}\xspace}
\newcommand{\xslabel}[2][]{\ifthenelse{\isempty{#1}}{\mathrm{#2}}{\mathrm{#2},#1}}
\newcommand{\Sigs}[1][]{\ensuremath{\Sigma_{\xslabel[#1]{s}}}\xspace}
\newcommand{\position}{\ensuremath{\vec{x}}\xspace}
\newcommand{\domg}{\dx[\Omega]}
\newcommand{\dx}[1][x]{\,d#1}
\newcommand{\Sigf}[1][]{\ensuremath{\Sigma_{\xslabel[#1]{f}}}\xspace}
\newcommand{\bnormal}{\ensuremath{\normal_\mathrm{b}}\xspace}
\newcommand{\normal}{\ensuremath{\vec{n}}\xspace}
\newcommand{\omgdnb}{\direction\cdot\bnormal}
\newcommand{\sfluxunit}{\,\ensuremath{\unit{cm}^{-2}\unit{s}^{-1}}\xspace}
\newcommand{\op}[1]{\mathbbl{#1}}
\newcommand{\voL}{{\op{\bold{L}}}}
\newcommand{\vect}[1]{\mathbf{#1}}
\newcommand{\oL}{{\op{{L}}}}
\newcommand{\voS}{{\op{\bold{S}}}}
\newcommand{\oS}{{\op{{S}}}}
\newcommand{\voF}{{\op{\bold{F}}}}
\newcommand{\oF}{{\op{{F}}}}

\newcommand*{\boundary}{\ensuremath{{\partial\domain}}\xspace}
\newcommand{\kernel}[1]{\mathbbm{#1}}

\newcommand{\ma}{\mathcal{A}}
\newcommand{\mb}{\mathcal{B}}
\newcommand{\mf}{\mathcal{F}}
\newcommand{\mj}{\mathcal{J}}

\newcommand{\norm}[2][{}]{\lVert#2\rVert_{#1}}
\newcommand{\vA}{\vect{A}}
\newcommand{\vM}{\vect{M}}
\newcommand{\vR}{\vect{R}}
\newcommand{\vP}{\vect{P}}
\newcommand{\vI}{\vect{I}}
\newcommand{\vx}{\vect{x}}
\newcommand{\vb}{\vect{b}}
\newcommand{\ve}{\vect{e}}
\newcommand{\vr}{\vect{r}}
\newcommand{\vz}{\vect{z}}
\newcommand{\DC}[1][]{\ensuremath{\mathrm{D}_{#1}}\xspace}

\newcommand{\Sigr}[1][]{\ensuremath{\Sigma_{\xslabel[#1]{r}}}\xspace}

\newcommand{\dvect}[1]{\vec{\mathbf{#1}}}

\begin{document}

\begin{frontmatter}

\title{Neutron transport criticality calculations using a parallel monolithic multilevel Schwarz 
preconditioner together with a  nonlinear diffusion acceleration method}

\author[firstauthor]{Fande Kong\corref{mycorrespondingauthor}}
\address[firstauthor]{Computational Frameworks, Idaho National Laboratory, P.O. Box 1625, Idaho Falls, ID 83415-3840}
\ead{fande.kong@inl.gov}
\cortext[mycorrespondingauthor]{Corresponding author}

\begin{abstract}
The multigroup neutron transport criticality calculations using modern supercomputers 
 have  been widely employed in a nuclear  reactor analysis  for 
 studying whether or not a system is self-sustaining.  
 However, the design and development of  efficient 
 parallel algorithms for the transport criticality calculations is challenging 
 especially when the number of processor   
 cores is large and an  unstructured  mesh is  adopted.  
In particular, both the compute time 
 and the memory usage  have to be carefully  taken 
 into consideration  due to the dimensionality  of the 
 neutron transport equations. In this paper, we study a monolithic multilevel 
 Schwarz preconditioner  for the transport criticality  calculations based on  
 a nonlinear diffusion acceleration (NDA) method. 
 In NDA, the multigroup  nonlinear diffusion equations are computed using 
 an inexact Jacobian-free Newton method with an initial guess generated from 
 a few inverse power iterations. The computed scalar fluxes and eigenvalue  are used 
 to evaluate  the fission and scattering terms of the transport equations, and then 
 the nonlinear system of  transport equations is simplified  to a linear system  of  equations. 
 The linear  systems of  equations arising from the discretizations  of the nonlinear 
 diffusion  equations and the transport equations need to be efficiently   solved.
We propose a monolithic  multilevel Schwarz method  that   
is capable of efficiently handling the  systems of linear equations 
 for both   the transport system  and  the diffusion system.
However, in the multilevel method,  
 algebraically constructing  coarse spaces is expensive and often unscalable.
 We study a  subspace-based coarsening 
 algorithm to  address such a challenge  
 by exploring  the matrix structures of  the transport equations 
 and the nonlinear diffusion equations. 
 We numerically demonstrate that the monolithic   multilevel preconditioner with 
 the subspace-based coarsening algorithm  is twice as fast as  that equipped 
 with  an unmodified  coarsening
 approach on  thousands of processor cores for 
 an unstructured mesh neutron transport problem with 
 billions of unknowns.

\end{abstract}

\begin{keyword}
Multilevel Schwarz preconditioner,  
nonlinear diffusion acceleration,  
multigroup neutron transport equations,  
parallel processing,
domain decomposition methods
\end{keyword}

\end{frontmatter}

\section{Introduction}
The neutron transport criticality calculations  play an important  role 
in a nuclear reactor analysis
\cite{lewis1984computational, duderstadt1976nuclear}. 
A neutron transport system is said to be critical if the rate of fission neutron production 
just equals to  the neutron losses because of absorption and leakage \cite{lewis1984computational}. 
A generalized eigenvalue problem  needs  to be calculated for checking the criticality   of 
the nuclear reactor system. The criticality   calculations of the neutron 
transport problems are challenging, and a tremendous amount of computational 
resources are  required
because of a high dimensionality; e.g., 1D energy, 2D angle and 3D spatial space. 
With an advancement of modern supercomputers,  the high-resolution simulations of the 
neutron transport problems become  possible 
\cite{kaushik2009enabling, kunen2015kripke, davidson2014massively}. 
However, for efficiently using  supercomputers, a scalable
parallel eigenvalue solver need to be designed and developed 
with taking both the compute time and the memory usage into consideration. 
In this work, we propose a highly parallel 
monolithic eigenvalue solver consisting  of  a nonlinear diffusion acceleration,
a Jacobian-free Newton-Krylov method and a multilevel 
 Schwarz preconditioner.

One of the simplest approaches for computing the criticality  of the multigroup  neutron 
transport problems  is inverse power iteration \cite{lewis1984computational, saad2011numerical},
and the inverse power iteration is still widely  used today due to its implementation simplicity. 
The inverse power iteration is often referred to as ``outer iteration". Sometimes, 
 the inverse power iteration is accelerated with 
DSA  (Diffusion Synthetic Acceleration) \cite{alcouffe1977diffusion}
or NDA (Nonlinear Diffusion Acceleration) \cite{smith2002full, schunert2017flexible}.
During  each inverse power  iteration, an inverse of the streaming-collision operator 
is required, and the calculation of the operator inverse can be generalized to solve  
a linear system of equations.  Traditionally, a Gauss-Seidel iterative method 
is employed for solving the linear system of equations group by group where
the scattering and the fission terms are computed using the previously   updated  solution 
\cite{lewis1984computational,kaushik2009enabling}.  For each Gauss-Seidel iteration,
a within-group subsystem of equations need to be efficiently solved.
Two of the most popular approaches for computing the within-group linear system 
of equations are transport sweeps  \cite{kunen2015kripke, davidson2014massively, shaner2016verification} 
and multilevel methods  
\cite{kong2019scalable, kong2019highly, kanschat2014robust, chang2007spatial, manteuffel1995fast}. 
In the current work, we consider multilevel preconditioned 
 iterative Krylov subspace methods that are easy to parallelize
 and can handle fully irregular meshes and different spatial discretizations.
The Krylov subspace methods mostly  involve  
matrix-vector multiplications  that can be efficiently implemented
 and optimized   for parallel calculations  \cite{petsc-user-ref, zhang2018vectorized}.
In addition, the multilevel   preconditioning   technique can be carried out in 
a scalable way based on the domain decomposition methods, and it is employed  to 
accelerate   the convergence 
of the  Krylov subspace methods, which will be discussed shortly.

However, 
the inverse power iteration based algorithm framework may converge slowly 
 when the smallest eigenvalue and the second smallest eigenvalue 
 are close to each other, which  
is often true for realistic problems. 
To address this challenge, while there exist other popular schemes 
such as  Wielandt shift \cite{saad2011numerical, yee2017space}, 
we employ  
an inexact Jacobian-free Newton
method \cite{dembo1982inexact, knoll2004jacobian, knoll2011acceleration}.
The Newton-based eigenvalue solver has been successfully  applied to 
many applications  
\cite{ kong2019scalable, kong2019highly, knoll2011acceleration, gill2011newton, calef2013nonlinear,  kong2018fully}, 
just name a few.
A direct application of the Jacobian-free Newton-based eigenvalue solver to the multigroup 
 neutron transport equations sometimes  could be expensive because the dimensionality of the nonlinear
 system of   
   equations is high.  One of  possible strategies  to improve the performance is to 
 apply a nonlinear diffusion acceleration (NDA) method  \cite{smith2002full, schunert2017flexible,wang2018rattlesnake}, 
 where
 the neutron transport problem is reformulated  as  a ``fixed source" problem by evaluating the scattering 
 and the fission terms using the scalar fluxes and eigenvalue of the diffusion equations. The NDA method 
 was originally proposed in \cite{smith2002full}, and its basic idea is directly related  to 
 the coarse mesh finite difference method investigated  in \cite{smith1983nodal}. The NDA method 
 combined with  multigrid in energy was recently studied in \cite{cornejo2016nonlinear, cornejo2019iteration}.

One of the most challenging tasks in the NDA method is 
to develop a robust and efficient preconditioner for both the transport system and the
diffusion system because the transport system is large 
and the diffusion system is nonsymmetric due to the closure terms. 
In this work, we study a  parallel monolithic multilevel preconditioner 
together with a subspace-based coarsening algorithm 
for both the transport system and the diffusion system. 
More precisely, an inexact Jacobian-free Newton method  
with an initial guess generated  from a few inverse power  iterations 
 is employed for the nonlinear diffusion system,  and during  each Newton iteration 
the Jacobian system is calculated using the monolithic multilevel  
Schwarz preconditioner together with 
GMRES, where unlike a Gauss-Seidel  type 
approach the monolithic  method handles  all energy  groups simultaneously. 
The transport system is reformed as a fixed source problem under 
the NDA context, and a linear system of transport equations 
is also solved using the monolithic multilevel preconditioner together with GMRES.
Note, there are other possibilities to explore the basic idea of multielvel methods;
e.g., the multilevel methods in energy were studied in  
\cite{cornejo2016nonlinear, cornejo2019iteration, slaybaugh2013multigrid};
the multilevel methods were  constructed in both spatial space and 
energy  \cite{yee2017multilevel, yee2019multilevel}.
We consider only  the algebraic spatial version of the multlevel methods in the current study.

 In the multilevel  methods, 
the construction of coarse spaces has  a critical  impact on the preconditioner performance. 
In general,  
a coarse space can be constructed either  geometrically
 \cite{kong2016highly,kong2017scalable} 
or algebraically \cite{kong2018fully}. 
We choose to use an algebraic coarsening algorithm in this work because 
it is general and has an out-of-box nature for end-users.
But it is well-known that the construction of coarse spaces is expensive in terms of the compute time and 
the memory usage  \cite{yang2002boomeramg},
and for some  applications the setup phase may be even more expensive than the solve phase. 
To overcome this difficulty, we introduce a subspace-based coarsening algorithm 
by exploring the structures  of the preconditioning matrices.
In the subspace-based coarsening algorithm, only a submatrix 
needs to be coarsened to generate coarse spaces, and therefore it is able 
to reduce the computational cost in time and in memory. 
The subspace-based coarsening was successfully applied to 
the unaccelerated neutron transport calculations in our previous 
work \cite{kong2019scalable, kong2019highly},
and it is improved and extended to solve both the nonlinear diffusion system 
and the accelerated   transport system in the current study.
We numerically demonstrate that the proposed algorithm 
is more efficient than the unmodified traditional multilevel methods  for solving 
both the accelerated  transport system and the diffusion system 
with  thousands of processor cores for a 3D unstructured  mesh problem
with billions of unknowns.

The rest of this paper is organized as follows. In Section 2, a 
nonlinear diffusion acceleration method 
 is presented to decouple the scattering 
term and the fission term  from the streaming-collision  operator in the 
neutron transport problem to speedup the simulations.  
A parallel monolithic 
multilevel Schwarz preconditioner together with a subspace-based 
coarsening algorithm  is 
described  in Section 3 for both the accelerated transport system and the diffusion system.  
In Section 4, some numerical tests 
are carefully studied to demonstrate  the performance of the 
monolithic  multilevel Schwarz preconditioner
equipped with the subspace-based coarsening algorithm.  
A few remarks and conclusions are drawn in Section 5.

\section{Nonlinear diffusion acceleration}
In this section, we present the multigroup neutron transport 
equations for the criticality  calculations, and then
study a nonlinear diffusion acceleration method where a system of diffusion equations is solved to 
speedup the transport simulations.

\subsection{Multigroup neutron transport equations}
 The multigroup neutron transport equations are 
used to describe the neutron interactions with the 
background materials  in a nuclear reactor system, 
and  they read as follows:
\begin{equation}\label{eq:eigenvalue}
\begin{array}{llll}
\displaystyle
 \direction\cdot\grad\Psi_g (\position, \direction)+
 \Sigt[g]\Psi_g (\position, \direction) & =   \displaystyle
 \sum_{g'=1}^G  \displaystyle \int_{\sphere} \Sigs[g' \rightarrow g]  f_{g' \rightarrow g}(\direction'\cdot\direction)\Psi_{g'} (\position, \direction')\domg'  \\
 & \displaystyle + \frac{1}{4\pi}\frac{\chi_g}{k}\sum_{g'=1}^G\nu\Sigf[g']\Phi_{g'}( \position), \text{ in } \domain \times \sphere, \\
\end{array}
\end{equation} 
where $\Psi_g$  $[\afluxunit]$ are the neutron angular fluxes defined on 
 $\domain \times \sphere$,
 $g=1,2,...,G$.  Here $G$ is the number of energy groups,  $\sphere$ is a unit sphere of neutron 
 flying directions and $\domain$ is a 3D spatial domain, 
 for example as shown in Fig.~\ref{fig:hierarch_partition_3D_80}. 
In~\eqref{eq:eigenvalue},  $\position \in \domain$  $\left[ \unit{cm} \right]$ 
denotes  the independent spatial variable,  and  $\direction \in \sphere$ 
represents the  independent angular variable.  $\Sigt[g]$ $[ \unit{cm}^{-1} ]$ 
is the macroscopic total cross section of the $g$th energy group, and 
 $\Sigs[g'\rightarrow g]$ $[ \unit{cm}^{-1} ]$ 
represents  the macroscopic scattering cross section from the $g'$th group
 to the $g$th group.  $f_{g'\rightarrow g}$  $[ \unit{st}^{-1} ]$ is the scattering phase 
function that redistributes  neutrons  from the incoming 
directions $\direction'$ to some  outgoing directions  $\direction$.
In the second term of the right hand side of~\eqref{eq:eigenvalue}, 
$\chi_g$ denotes the spectrum for both prompt and delayed,   
$\nu$ is the averaged number of neutrons  emitted per fission,
$\Sigf[g]$  $[ \unit{cm}^{-1} ]$ represents  the macroscopic fission cross section, and
$\Phi_{g}$ $[\sfluxunit ]$ are the scalar fluxes  that can be calculated by taking an integral 
of   the angular fluxes over the unit  sphere
$$
\Phi_{g}  \equiv \int_\sphere \Psi_g  \domg.
$$

$k$ is  the inverse of  eigenvalue 
(the largest $k$ is  referred to as the neutron multiplication factor).  The neutron
transport criticality calculations focus  on   the largest value of $k$ and
 the corresponding eigenvalue vector, referred to as the fundamental mode.
The first term of~\eqref{eq:eigenvalue} is  the \emph{streaming term}, 
the second is the \emph{collision} term,
the third  is  the \emph{scattering term} and the forth  is the \emph{fission} term.  The scattering term 
couples  energy groups  through the scattering matrix, and connects  angular directions through the integral. 
Similarly, the fission term couples the angular directions together.  
The angular  and  energy 
coupling makes the neutron transport  criticality calculations  challenging. 

To rewrite~\eqref{eq:eigenvalue} as a  vector form,  
we introduce some notations  as follows:
\begin{subequations}
\begin{equation}\label{define:angular}
 \vect{\Psi} \equiv
 \begin{bmatrix}
 \displaystyle
 \Psi_1, 
  \displaystyle
 \Psi_2,
   \displaystyle
   \cdots,
   \displaystyle
 \Psi_G
\end{bmatrix}^T, 
\end{equation}
\begin{equation}\label{define:Loperator}
\voL \vect{\Psi} \equiv
 \begin{bmatrix}
 \displaystyle
  \oL_1 \Psi_1, 
  \displaystyle
   \oL_2 \Psi_2,
   \displaystyle
   \cdots,
   \displaystyle
  \oL_G \Psi_G
\end{bmatrix}^T, ~ \oL_g \Psi_g \equiv   \grad \cdot \direction \Psi_g+ \Sigt[g] \Psi_g, 
\end{equation}
\begin{equation}\label{define:Soperator}
\voS \vect{\Psi} \equiv
 \begin{bmatrix}
\displaystyle
 \oS_1 {\Psi_1},
  \displaystyle
  \oS_2 {\Psi_2},
  \displaystyle
  \cdots,
  \displaystyle
\oS_G {\Psi_G}
 \end{bmatrix}^T,  ~ \oS_g {\Psi_g} \equiv  \sum_{g'=1}^G\int_\sphere \Sigs[g'\rightarrow g]f_{g'\rightarrow g}\Psi_{g'}\domg', 
\end{equation}
\begin{equation}\label{define:Foperator}
\voF \vect{\Psi}  \equiv
 \begin{bmatrix}
 \displaystyle
\oF_1  \Psi_1 , 
  \displaystyle
\oF_2  \Psi_2  ,
  \cdots,
  \displaystyle
\oF_G  \Psi_G 
 \end{bmatrix}^T,  ~\oF_g  \Psi_g \equiv  \frac{1}{4\pi}\chi_{g}\sum_{g'=1}^G\nu\Sigf[g']\Phi_{g'}.
 \end{equation}
\end{subequations}
Here $\vect{\Psi}$ is the angular flux vector consisting of all energy groups.  $\voL$, 
$\voS$ and $\voF$ are the streaming-collision operator, 
the scattering operator and the fission  operator,
respectively.  With the operator notations, we rewrite~\eqref{eq:eigenvalue}  as
\begin{equation}\label{eq:operators}
\voL \vect{\Psi} = \voS \vect{\Psi} + \frac{1}{k}\voF \vect{\Psi}.
\end{equation}

Following a standard finite element technique, 
the weak form of~\eqref{eq:operators}  reads as
\begin{equation}\label{eq:unstabilized-bilinear}
\left(\voL^\ast\vect{\Psi}^\ast, \vect{\Psi}\right)
 +\left\langle \vect{\Psi}^\ast, \bar{\vect{\Psi}} \right\rangle^+
 -\left\langle \vect{\Psi}^\ast, \bar{\vect{\Psi}} \right\rangle^- =
 \left(\vect{\Psi}^\ast, \voS \vect{\Psi}\right)
  + \frac{1}{k}  \left (\vect{\Psi}^\ast , \voF \vect{\Psi} \right ), 
\end{equation}
where $\vect{\Psi}^\ast$ is the test function, and $\voL^\ast$ is the adjoint operator of $\voL$ defined 
as 
$$
\voL^\ast \vect{\Psi} \equiv
 \begin{bmatrix}
 \displaystyle
  \oL_1^\ast \Psi_1, 
  \displaystyle
   \oL_2^\ast \Psi_2,
   \displaystyle
   \cdots,
   \displaystyle
  \oL_G^\ast \Psi_G
\end{bmatrix}^T, ~ \oL_g^\ast \Psi_g \equiv  -  \grad \cdot \direction \Psi_g+ \Sigt[g] \Psi_g. 
$$
$(\cdot, \cdot)$ denotes a  function inner product over $\domain \times \sphere$,  for example,
$$
 \left( \vect{\Psi}^\ast, \vect{\Psi} \right) \equiv \sum_{g=1}^G\int_{\sphere}\domg\int_{\domain}\dx\, \Psi_g^\ast(\position,\direction)\Psi_g(\position, \direction).
$$
$\left\langle \cdot,  \cdot \right\rangle$ represents the inner product over the boundaries, that is,
$$
\begin{array}{llll}
\displaystyle   \left\langle  \vect{\Psi}^\ast, \bar{\vect{\Psi}} \right\rangle^\pm \equiv 
\sum_{g=1}^G \oint_{\boundary}\dx\int_{\sphere^\pm_{\bnormal}}\domg\, \left|\omgdnb\right| \Psi_g^\ast(\position,\direction)\bar{\Psi}_g(\position, \direction). 
\end{array}
$$
Here $\partial \domain$ is the boundary of  $\domain$,    $\bnormal$ is the outward unit normal vector on $\partial \domain$, and $\sphere^\pm_{\bnormal}=\{\direction \in \sphere :   \direction\cdot \bnormal  \gtrless 0 \}.$ 
The vacuum and reflecting boundary conditions are employed in this paper. 
For the vacuum boundary condition, $\bar{\vect{\Psi}}$ is expressed as 
 $$
 \bar{\vect{\Psi}} =
 \left \{
 \begin{array}{lll}
{\vect{\Psi}},  &\text{~on~}  \partial \domain, ~\direction\cdot \bnormal \geq 0,\\
 \vect{0}, & \text{~on~}  \partial \domain, ~\direction\cdot \bnormal <0.
 \end{array}
 \right .
$$
For the reflecting boundary condition,   $\bar{\vect{\Psi}}$ is  defined as 
 $$
 \bar{\vect{\Psi}} =
 \left \{
 \begin{array}{lll}
{\vect{\Psi}},  &\text{~on~}  \partial \domain, ~\direction\cdot \bnormal \geq 0,\\
 \vect{\vect{\Psi}}_r, & \text{~on~}  \partial \domain, ~\direction\cdot \bnormal <0,
 \end{array}
 \right .
$$
where $\vect{\vect{\Psi}}_r$ is the reflecting 
angular fluxes of  ${\vect{\Psi}}$ on  $\direction_r = \direction - 2(\omgdnb)\bnormal$.
The weak form~\eqref{eq:unstabilized-bilinear} is unstable, and certain stabilizing  techniques 
are  required. A stabilizing technique  called SAAF (self-adjoint angular flux)
 \cite{morel1999self, liscum2002finite} is chosen 
in this work, and  the SAAF is very similar to the well-known method 
SUPG (Streamline upwind/Petrov-Galerkin) \cite{brooks1982streamline}.
We ignore the description  of  SAAF here since it has been presented 
in other literatures \cite{kong2019scalable,kong2019highly,wang2018rattlesnake, morel1999self, liscum2002finite}. 

Without introducing any confusion,
the same weak  form notation  is also reused to represent its stabilized version. 
To simplify the description,
the weak form is rewritten as
\begin{equation}\label{eq:unstabilized-bilinear-simple}
\kernel{l} (\vect{\Psi}^\ast, \vect{\Psi}) - \kernel{s}_{\text{h}} (\vect{\Psi}^\ast, \vect{\Psi})  =
\kernel{s}_{\text{0}} (\vect{\Psi}^\ast, \vect{\Phi}) + \frac{1}{k} \kernel{f}  (\vect{\Psi}^\ast, \vect{\Phi}),
\end{equation}
with 
$$
\kernel{l} (\vect{\Psi}^\ast, \vect{\Psi}) \equiv  \left(\voL^\ast\vect{\Psi}^\ast, \vect{\Psi}\right)
 +\left\langle \vect{\Psi}^\ast, \bar{\vect{\Psi}} \right\rangle^+ -\left\langle \vect{\Psi}^\ast, \bar{\vect{\Psi}} \right\rangle^-,
$$
$$
\kernel{f}  (\vect{\Psi}^\ast, \vect{\Phi}) \equiv \left (\vect{\Psi}^\ast , \voF \vect{\Psi} \right ),
$$
and 
$$
\kernel{s}_{\text{h}} (\vect{\Psi}^\ast, \vect{\Psi}) +\kernel{s}_{\text{0}} (\vect{\Psi}^\ast, \vect{\Phi}) 
\equiv  \left(\vect{\Psi}^\ast, \voS \vect{\Psi}\right).
$$
Here $\kernel{s}_{\text{h}}$ corresponds to the high order moment fluxes, and 
$\kernel{s}_{\text{0}}$ corresponds to the zero  order moment fluxes  
that are completely determined 
by the scalar fluxes. 

\subsection{Diffusion-based acceleration technique}
\eqref{eq:unstabilized-bilinear-simple}  involves the scalar fluxes  
in the fission term and the scattering term.
An alternative way to compute the scalar fluxes 
is to solve  the multigroup neutron diffusion equations.
The resulting scalar fluxes and eigenvalue  are  
used to evaluate the fission term and the scattering term of
the neutron transport equations.  
In such a way, the multigroup neutron transport equations are simplified  to 
a linear problem  from an  eigenvalue problem. 
It is much cheaper to solve the linear system of 
equations  instead of the eigenvalue problem. 

The scalar fluxes are governed by the multigroup 
neutron diffusion equations,
\begin{equation} \label{eq:diffusion}
 - \grad\cdot\DC[g]\grad\Phi_g (\position) + \Sigr[g]\Phi_g  (\position)
 =  \sum_{g'\neq g}\Sigs[g'\rightarrow g]\Phi_{g'}(\position) 
 + \frac{\chi_g}{k}\sum_{g'=1}^G\nu\Sigf[g']\Phi_{g'} (\position), {\text{~in~}} \domain,
\end{equation}
where  $\DC[g] [\unit{cm}] $ is the   diffusion coefficient, 
and  $\Sigr[g] \equiv \Sigt[g] - \Sigs[g\rightarrow g]$, in $[\unit{cm}^{-1}]$, 
is the macroscopic removal cross section.  The first  term 
is the \emph{diffusion term} and the second  is the \emph{removal term}.
The first  term on the right hand side is the \emph{scattering term}, 
which couples the scalar fluxes  together through the scattering matrix.
The last term is the \emph{fission term}.
For conciseness, we also write Eq.~\eqref{eq:diffusion} in a vector form,
\begin{equation}\label{eq:vecdiffusion}
- \grad\cdot\dvect{J} + \voL_{\text{r}}\vect{\Phi} = \voS_\text{d}\vect{\Phi} +  \frac{1}{k} \voF_0 \vect{\Phi},
\end{equation}
where 
$$
\dvect{J} \equiv [ \vec{J}_1,  \vec{J}_2, ...,  \vec{J}_G]^T, ~ \vec{J}_g \equiv \DC[g]\grad\Phi_g,
$$
$$
\voL_r \vect{\Phi} \equiv [\op{L}_{\text{r},1}  \Phi_{1},  \op{L}_{\text{r},2}  \Phi_{2},..., \op{L}_{\text{r},G}  \Phi_{G}] ^T, 
\op{L}_{\text{r},g}  \Phi_{g} \equiv   \Sigr[g] \Phi_{g},
$$
$$
\voS_\text{d}\vect{\Phi} \equiv [\op{S}_{\text{d},1}  \Phi_{1},  \op{S}_{\text{d},2}  \Phi_{2},..., \op{S}_{\text{d},G}  \Phi_{G} ]^T, 
 \op{S}_{\text{d},g} \equiv  \sum_{g'\ne g} \Sigs[g'\rightarrow g] \Phi_{g'},
$$
and
$$
\voF_0 \vect{\Phi}  \equiv
 \begin{bmatrix}
 \displaystyle
\oF_{0,1}  \Phi_1 , 
  \displaystyle
\oF_{0,2}  \Phi_2  ,
  \cdots,
  \displaystyle
\oF_{0,G}  \Phi_G 
 \end{bmatrix}^T,  ~\oF_{0,g}  \Phi_g \equiv \chi_{g}\sum_{g'=1}^G\nu\Sigf[g']\Phi_{g'}.
 $$
Multiply~\eqref{eq:vecdiffusion} by a test function $\vect{\Phi}^\ast$ and integrate  
by parts, then we have the weak form of~\eqref{eq:vecdiffusion}
as follows
\begin{equation}\label{weakform:diffusion}
 \kernel{b}_\text{diff} (\vect{\Phi}^\ast, \vect{\Phi}) =   \frac{1}{k}  \kernel{f}_\text{diff} (\vect{\Phi}^\ast, \vect{\Phi}),
\end{equation}
with 
$$
   \kernel{b}_\text{diff}(\vect{\Phi}^\ast, \vect{\Phi})  =  
  \left(\grad\vect{\Phi}^\ast, \op{D}\grad\vect{\Phi}\right)_\domain +
  \left<\vect{\Phi}^\ast, \frac{1}{4}\vect{\Phi}\right>_{\boundary} +
  \left(\vect{\Phi}^\ast, \voL_{\text{r}}\vect{\Phi}\right)_\domain 
  \displaystyle - 
  \left(\vect{\Phi}^\ast, \voS_\text{d}\vect{\Phi}\right)_\domain,
 $$
 and 
 $$
\kernel{f}_\text{diff} (\vect{\Phi}^\ast, \vect{\Phi}) =   \left(\vect{\Phi}^\ast, \voF_0\vect{\Phi}\right)_\domain, 
$$
where  $\left(\cdot,  \cdot \right)_\domain$  denotes a function inner product over $\domain$, 
and $ \left<\cdot, \cdot \right >_{\boundary}$ represents the function inner product on $\boundary$.
In~\eqref{weakform:diffusion}, the Robin boundary condition is used,
$$
 \frac{\Phi_g}{4}+\DC[g]\grad\Phi_g\cdot\bnormal = 0,  \text{~on~}\boundary.
$$

For accelerating the solution of~\eqref{eq:unstabilized-bilinear-simple},
the low-order diffusion equations~\eqref{weakform:diffusion} need additional nonlinear 
closure terms so that the low-order scalar fluxes become identical to 
the projection of the high order solution into the low order space when the algorithm
converges.  The following closure terms are  employed when~\eqref{eq:unstabilized-bilinear-simple}
  is stabilized with the SAAF method and  discretized in angle  using the discrete ordinates (SN),
\begin{equation}\label{eq:closureterm}
\kernel{C}[\vect{\Psi}](\vect{\Phi}^\ast, \vect{\Phi}) \equiv \left(\grad\vect{\Phi}^\ast, \tilde{\dvect{D}}\vect{\Phi} \right)_\domain + \left<\vect{\Phi}^\ast, \boldsymbol{\gamma}\vect{\Phi} \right>_{\partial\domain}.
\end{equation}
Here 
$$
\tilde{\dvect{D}} \equiv [\tilde{\vec{D}}_1, \tilde{\vec{D}}_2,..., \tilde{\vec{D}}_G]^T,
$$
with
$$
\footnotesize
\displaystyle  \tilde{\vec{D}}_g = \displaystyle   \frac{ \displaystyle    \int_\sphere\left( \tau_g\direction\direction\cdot\grad\Psi_g + (\tau_g\Sigt[g] - 1)\direction\Psi_g - \tau_g\sum_{g'=1}^G \Sigs[1]^{g'\rightarrow g}\direction\Psi_{g'}  - \DC[g]\grad\Psi_g\right)\domg}{ \displaystyle   \int_\sphere \Psi_g\domg}, \\
$$
and 
$$
\displaystyle \boldsymbol{\gamma} \equiv [ \gamma_1,  \gamma_2, ...,  \gamma_G]^T,
$$
with
$$
\displaystyle  \gamma_g = \frac{\displaystyle \int_{\omgdnb>0} \left|\omgdnb\right| \Psi_g\domg}{\displaystyle \int_\sphere \Psi_g\domg} - \frac{1}{4}.
$$
Here $\Sigs[1]^{g'\rightarrow g}$ is the first order of scattering 
cross section (see more details in \cite{wang2018rattlesnake}),  
and $\tau_g$ is the stabilization parameter defined as 
\begin{equation} \label{eq:tau}
 \tau_g = \left\{\begin{array}{lr}\frac{1}{c\Sigt[g]}, & ch\Sigt[g]\geq \varsigma, \\ \frac{h}{\varsigma}, & ch\Sigt[g] < \varsigma,\end{array}\right.
\end{equation}
where $h$ is the characteristic  length of a mesh  element,  $\varsigma$ is usually chosen to be a constant of  0.5 and 
$c$ is a constant as well (1.0 by default).  A detailed explanation  of \eqref{eq:tau} can be found in \cite{wang2018rattlesnake}.
With  these notations,  the nonlinear diffusion acceleration 
scheme is summarized in Alg.~\ref{alg:nda}.
\begin{algorithm}
  \caption{Nonlinear diffusion acceleration method
    \label{alg:nda}}
  \begin{algorithmic}[1]
    \State Solve the low-order diffusion system  without the closure terms:
    $$
     \kernel{b}_\text{diff} (\vect{\Phi}^\ast, \vect{\Phi}^{(0)}) =   \frac{1}{k^{(0)}}  \kernel{f}_\text{diff} (\vect{\Phi}^\ast, \vect{\Phi}^{(0)})
   $$
    \State $n= 0$
    \For { $n<\max_{\text{nda}}$ and $\epsilon_1 > tol$}
      \State Solve  the transport system:
\begin{equation}\label{eq:nda-transport}
               \kernel{l} (\vect{\Psi}^\ast, \vect{\Psi}^{(n+1)}) - \kernel{s}_{\text{h}} (\vect{\Psi}^\ast, \vect{\Psi}^{(n)})  =
               \kernel{s}_{\text{0}} (\vect{\Psi}^\ast, \vect{\Phi}^{(n)}) + \frac{1}{k^{(n)}} \kernel{f}  (\vect{\Psi}^\ast, \vect{\Phi}^{(n)})
\end{equation}
\State  Update the coefficients, $\tilde{\dvect{D}}$, of the closure terms
      \State Solve the diffusion  system with the closure terms:
\begin{equation}\label{eq:nda-diffusion}
         \kernel{b}_\text{diff}(\vect{\Phi}^\ast, \vect{\Phi}^{(n+1)}) + \kernel{C}[\vect{\Psi}^{{(n+1)}}](\vect{\Phi}^\ast, \vect{\Phi}^{(n+1)}) = \frac{1}{k^{(n+1)}}\kernel{f}_\text{diff}(\vect{\Phi}^\ast, \vect{\Phi}^{(n+1)})
\end{equation}
      \State Calculate  $\epsilon_1 = \left| \vect{\Phi}^{(n+1)} - \vect{\Phi}^{(n)} \right|$
      \State $n += 1$
   \EndFor
  \end{algorithmic}
\end{algorithm}

The low order diffusion equations~\eqref{eq:nda-diffusion} are discretized  
using the first order continuous  finite element  \cite{kong2018fully, wang2018rattlesnake}, 
and the multigroup transport equations~\eqref{eq:nda-transport} are discretized 
with the first order finite element 
method in space  and with the SN scheme   
in angle \cite{kong2019scalable,kong2019highly, wang2018rattlesnake}.
After the spatial and angular discretizations, the multigroup transport equations at 
line 4 of Alg.~\ref{alg:nda} correspond to 
a linear system of equations since the fission term and part of the scattering term 
are evaluated using the scalar fluxes and the eigenvalue  computed in  the low order 
diffusion system. 
An algebraic  generalized  eigenvalue system is generated after the spatial discretization 
 of the low order diffusion system. 
The low order diffusion  system is challenging to solve 
since the extra closure terms can have 
a significant impact on diffusion coefficients,
and meanwhile 
the high order transport  system is expensive to solve 
due to the high dimensionality. 
We next study a novel monolithic multilevel  Schwarz 
method that is capable of solving both 
the low order diffusion system and  the high order transport system 
while constructing coarse spaces in an efficient way. 

\section{Monolithic  multilevel Schwarz preconditioner}
In this section, we briefly describe a parallel  monolithic  
 algorithm framework, and then 
study a   multilevel Schwarz preconditioner
equipped with a subspace-based 
coarsening approach for both the transport equations 
and the diffusion equations. 

\subsection{Parallel Newton accelerated eigenvalue solver}
After the spatial discretization of~\eqref{eq:nda-diffusion},  
an algebraic generalized eigenvalue system 
reads as 
\begin{equation}\label{eq:algebraic-diffusion}
\ma \vect{\Phi}  = \frac{1}{k}\mb\vect{\Phi},
\end{equation}
where $\ma$ is the matrix consisting of the first and the second 
terms of~\eqref{eq:nda-diffusion}, and $\mb$ represents the matrix 
corresponding to the fission term. $\ma$ and $\mb$ are not formed 
explicitly for saving memory. The notation $\vect{\Phi}$ is reused to denote 
the discretized  version of the scalar fluxes, that is,  
it is a vector whose components are 
the finite element nodal values. The neutron transport criticality  calculations 
target at  the largest value of $k$ in~\eqref{eq:algebraic-diffusion}. 

One of the simplest approaches  to calculate  
the largest value of $k$
  is inverse power iteration \cite{lewis1984computational, saad2011numerical}, 
which starts with an initial guess and computes  a new  approximation 
with solving a linear system of equations.  More precisely, with a given initial pair 
  $( \vect{\Phi}_0,  k_0 =  \norm{\mb  \vect{\Phi}_{0}})$,  
  a new approximation $(\vect{\Phi}_{n+1}, k_{n+1})$, $n=0,1,..., \max_{\text{eigen}} $, 
  is obtained using  the 
  following steps:
  \begin{subequations} \label{alg:inversepower}
\begin{equation}\label{eq:inversepower1}
\ma \vect{\Phi}_{n+1} = \mb \vect{\Phi}_{n},
\end{equation}
\begin{equation}\label{eq:inversepower2}
k_{n+1}  = \norm{\mb  \vect{\Phi}_{n+1}},
\end{equation}
\begin{equation}\label{eq:inversepower3}
\displaystyle \mb \vect{\Phi}_{n+1}  \leftarrow \frac{1} {k_{n+1}} \mb \vect{\Phi}_{n+1}.
\end{equation}
\end{subequations}
Here the first step solves a linear system of equations,  
the second step calculates  the norm of the 
right hand size as a new eigenvalue approximation, 
and finally 
the right hand side is scaled in place. 
 
The inverse power iteration converges slowly when the smallest 
and the second smallest eigenvalues are close to each other,
which often occurs in realistic problems. To overcome this 
difficulty,  an inexact Newton method \cite{dembo1982inexact} is employed 
to accelerate the inverse power iteration by 
reforming~\eqref{alg:inversepower} as a nonlinear system of equations,
\begin{equation}\label{nonlinear_equation_newton}
 \mf(\vect{\Phi}) \equiv
\begin{array}{llll}
\displaystyle  \ma \vect{\Phi} - \frac{1}{\norm{\mb \vect{\Phi}}} \mb \vect{\Phi}=0. \\
\end{array}
\end{equation}
Newton  starts with an initial guess $\vect{\Phi}_0$ that, in this paper, is calculated with a few inverse power iterations,
and then a new   approximation solution,   $\vect{\Phi}_{n+1}$, $n=0,1,..,\max_{\text{newton}}$,  is updated   as follows,
\begin{equation}\label{alg:newton_update}
\vect{\Phi}_{n+1} =  \vect{\Phi}_{n} + \alpha_n  \Delta {\vect{\Phi}}_{n}.
\end{equation}
Here  $\alpha_n$ is a Newton step length  calculated using a backtracking line search scheme \cite{dennis1996numerical},  and 
$\Delta {\vect{\Phi}}_{n}$ is obtained by solving the Jacobian system
\begin{equation}\label{alg:newton_jacobian}
\mj( \vect{\Phi}_{n}) \Delta{\vect{\Phi}}_{n}  = -\mf(\vect{\Phi}_{n}),
\end{equation}
where $\mj( \vect{\Phi}_{n})$ is the Jacobian matrix evaluated at $\vect{\Phi}_{n}$ and 
$\mf(\vect{\Phi}_{n})$ is the nonlinear function residual at $\vect{\Phi}_{n}$. 
Explicitly forming $\mj$ is memory intensive since 
there exists  significant coupling among 
 variables through the scattering term and the fission term. 
  In addition,  the derivatives of $1/\norm{\mb \vect{\Phi}}$ is difficult  
to compute. To fix these issues, $\mj$  is carried out in a matrix-free manner,
that is, the action of $\mj$ on a vector $ \Delta{\vect{\Phi}}_{n}$ is 
implemented via finite difference
$$
\mj(\vect{\Phi}_{n}) \Delta{\vect{\Phi}}_{n} = \frac{\mf(\vect{\Phi}_{n}+ \beta  \Delta{\vect{\Phi}}_{n}  ) - \mf(\vect{\Phi}_{n})}{\beta},
$$
where $\beta$ is a small   parameter calculated using a strategy  proposed in \cite{pernice1998nitsol}.   
The resulting inexat Newton  is 
often referred to as ``Jacobian-free Newton method" \cite{knoll2004jacobian, knoll2011acceleration}
that is also very useful when the Jacobain is difficult to form due to 
for example the complicated physics  procedures  under a  multiphysics environment. We summarize 
the  Newton accelerated eigenvalue solver in Alg.~\ref{alg:newton-eigenvalue-solver}, 
and will discuss the preconditioning technique    
in next section.  
 \begin{algorithm}
  \caption{Newton accelerated eigenvalue solver
    \label{alg:newton-eigenvalue-solver}}
  \begin{algorithmic}[1]
   \State Initialize:   $\tilde{\vect{\Phi}}_0, ~\tilde{k}_0 =  \norm{\mb \tilde{ \vect{\Phi}}_{0}}$
   \For {$n=0,1,2,..., \max_{\text{eigen}}$}
       \State Solve a linear sytem: $\ma \tilde{\vect{\Phi}}_{n+1} = \mb \tilde{\vect{\Phi}}_{n}$
       \State Update eigenvalue: $\tilde{k}_{n+1}  = \norm{\mb  \tilde{\vect{\Phi}}_{n+1}}$
       \State  Scale right hand side:  $ \mb \tilde{\vect{\Phi}}_{n+1}  \leftarrow \frac{1} {k_{n+1}} \mb \tilde{\vect{\Phi}}_{n+1}$
   \EndFor
   \State Set an initial guess  for Newton:  ${\vect{\Phi}}_0 = \tilde{\vect{\Phi}}_{n}$
   \State n=0
   \For{$n<\max_{\text{newton}}$ and $\epsilon > \text{tol}_{\text{newton}}$}
   \State  Solve the Jacobian sytem: $\mj( \vect{\Phi}_{n}) \Delta{\vect{\Phi}}_{n}  = -\mf(\vect{\Phi}_{n})$
   \State  Perform a line search to compute $\alpha_n$
   \State  Update solution: $\vect{\Phi}_{n+1} =  \vect{\Phi}_{n} + \alpha_n  \Delta {\vect{\Phi}}_{n}$
   \State  $\epsilon = \norm{\mf(\vect{\Phi}_{n+1}}/\norm{\mf(\vect{\Phi}_{0}}$
   \State   n += 1
   \EndFor
  \end{algorithmic}
\end{algorithm}
Note that the Newton-Krylov method  has been widely and successfully  employed for nonlinear systems 
of equations arising from different engineering  areas; e.g., elasticity  problems \cite{kong2016highly},
incompressible flows \cite{kong2018efficient}, 
fluid-structure interactions \cite{kong2017scalable, kong2018scalability, kong2016parallel, kong2019simulation}.  
In this work, the Newton-Krylov method is employed to solve the eigenvalue problems.

\subsection{Single-component grid based coarse spaces}
In this section, we present a monolithic  multilevel method for~\eqref{eq:nda-transport}, 
\eqref{eq:inversepower1} and~\eqref{alg:newton_jacobian},
where a subspace-based coarsening algorithm is studied to 
generate a sequence of  Single-component Grid based coarse spaces
(SG coarse spaces) for 
improving the parallel   performance.

Three linear systems of equations, \eqref{eq:nda-transport}, 
\eqref{eq:inversepower1} and~\eqref{alg:newton_jacobian}, 
need to  be solved, and we hence denote them 
in an unified general linear system of equations,
\begin{equation}\label{eq:linear-system}
\vA\vx = \vb,
\end{equation}
where $\vA = \ma$ in~\eqref{eq:inversepower1},  
$\vA=\mj$ in~\eqref{alg:newton_jacobian}, 
and $\vA$ corresponds to  the discretized  version of $\kernel{l}$ 
in~\eqref{eq:nda-transport}.  $\vx$ is the nodal values of $\vect{\Phi}$ 
in~\eqref{eq:inversepower1} and~\eqref{alg:newton_jacobian}, 
and corresponds to the discretized 
version of $\vect{\Psi}$ in~\eqref{eq:nda-transport}.
The linear system~\eqref{eq:linear-system} is  
solved using an iterative Krylov subspace method; 
e.g., GMRES \cite{saad1986gmres}, and  
an efficient   preconditioner  is required to speedup the convergence.  
In the current study, a standard restarted 
GMRES from PETSc \cite{petsc-user-ref} is chosen, 
and no change was made to GMRES.
The restart number of GMRES is 30 by default in PETSc.
While other Krylov subspace methods \cite{saad2003iterative} are possible
for solving the Jacobian system, we consider only GMRES  in the 
current study. GMRES was chosen as the default linear solver 
in the multiphysics framework MOOSE \cite{permann2019moose} because it 
appears sufficiently general for a wide range of multiphysics problems. 

The right preconditioned linear system is rewritten as 
\begin{equation}\label{eq:preconditioned-linear-system}
\vA \vM^{-1} \vM  \vx = \vb,
\end{equation}
where $\vM$ is a preconditioning matrix, and $\vM^{-1}$ is a preconditioning process.  
$\vM$ is an approximation to 
$\ma$ in~\eqref{eq:inversepower1},  to $\mj$ in~\eqref{alg:newton_jacobian},
 and $\vM$ corresponds to the discretized version of $\kernel{l}$ in~\eqref{eq:nda-transport}.
The preconditioned system~\eqref{eq:preconditioned-linear-system} 
is computed in  two steps,
\begin{subequations}\label{eq:two-step-preconditioner}
\begin{equation}\label{eq:two-step-preconditioner1}
\vA \vM^{-1} \tilde{\vx} = \vb,
\end{equation}
\begin{equation}\label{eq:two-step-preconditioner2}
 \vM  \vx =  \tilde{\vx},
\end{equation}
\end{subequations}
where $\tilde{\vx} $ is an intermediate auxiliary vector. 
In~\eqref{eq:two-step-preconditioner}, the preconditioner is 
applied from the right side, but we want to mention that the preconditioner 
can be applied from the left side as well. In general,  
a preconditioning  process aims at a correction with solving 
the residual equations,
\begin{equation}\label{alg:newton_preconditioner_process}
\vM \ve = \vr,
\end{equation}
where $\ve$ is a correction vector, and $\vr$ is the residual vector from the outer solver, GMRES. 

\begin{figure}
 \centering
 \includegraphics[width=0.49\linewidth]{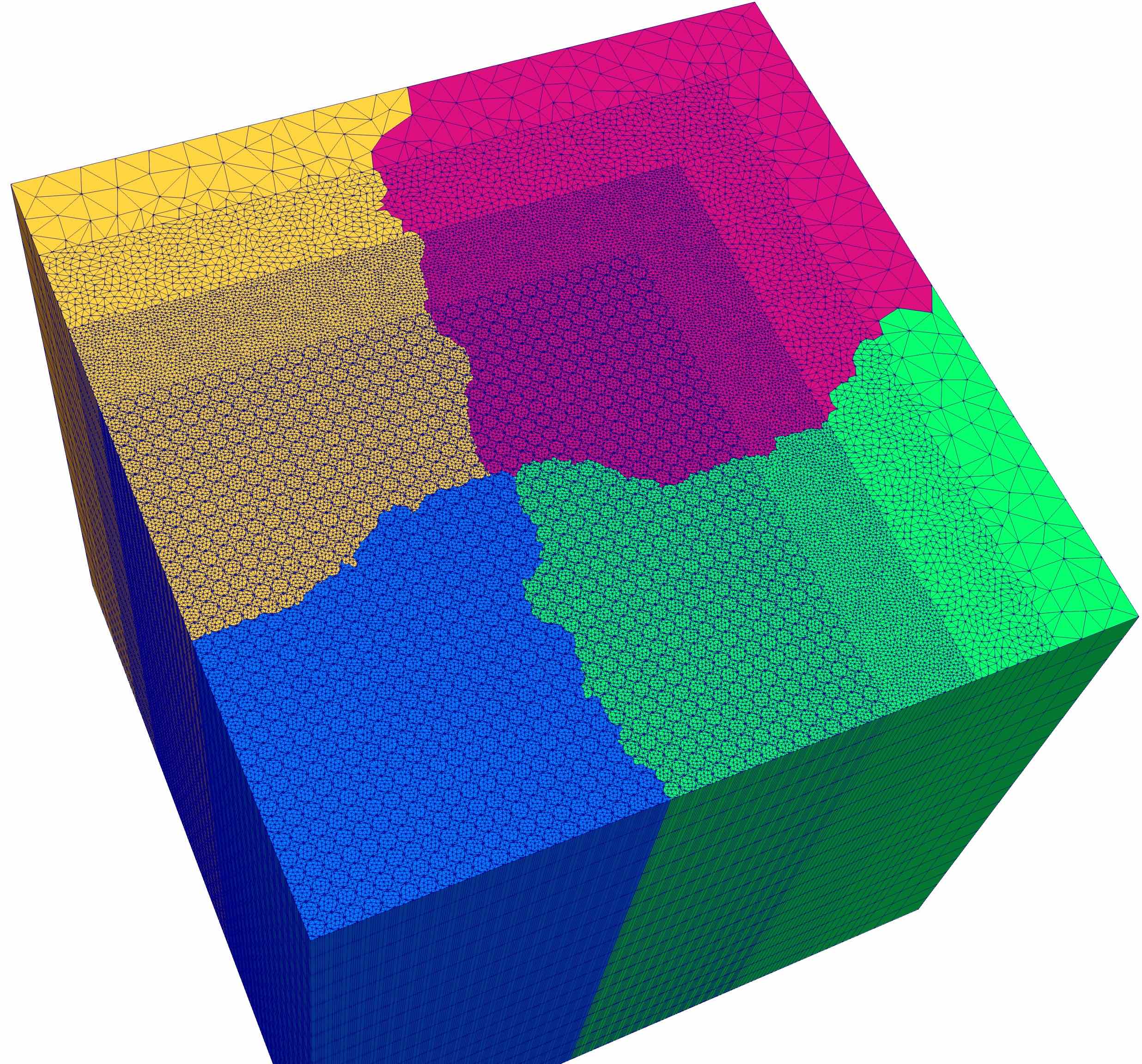}
 \includegraphics[width=0.49\linewidth]{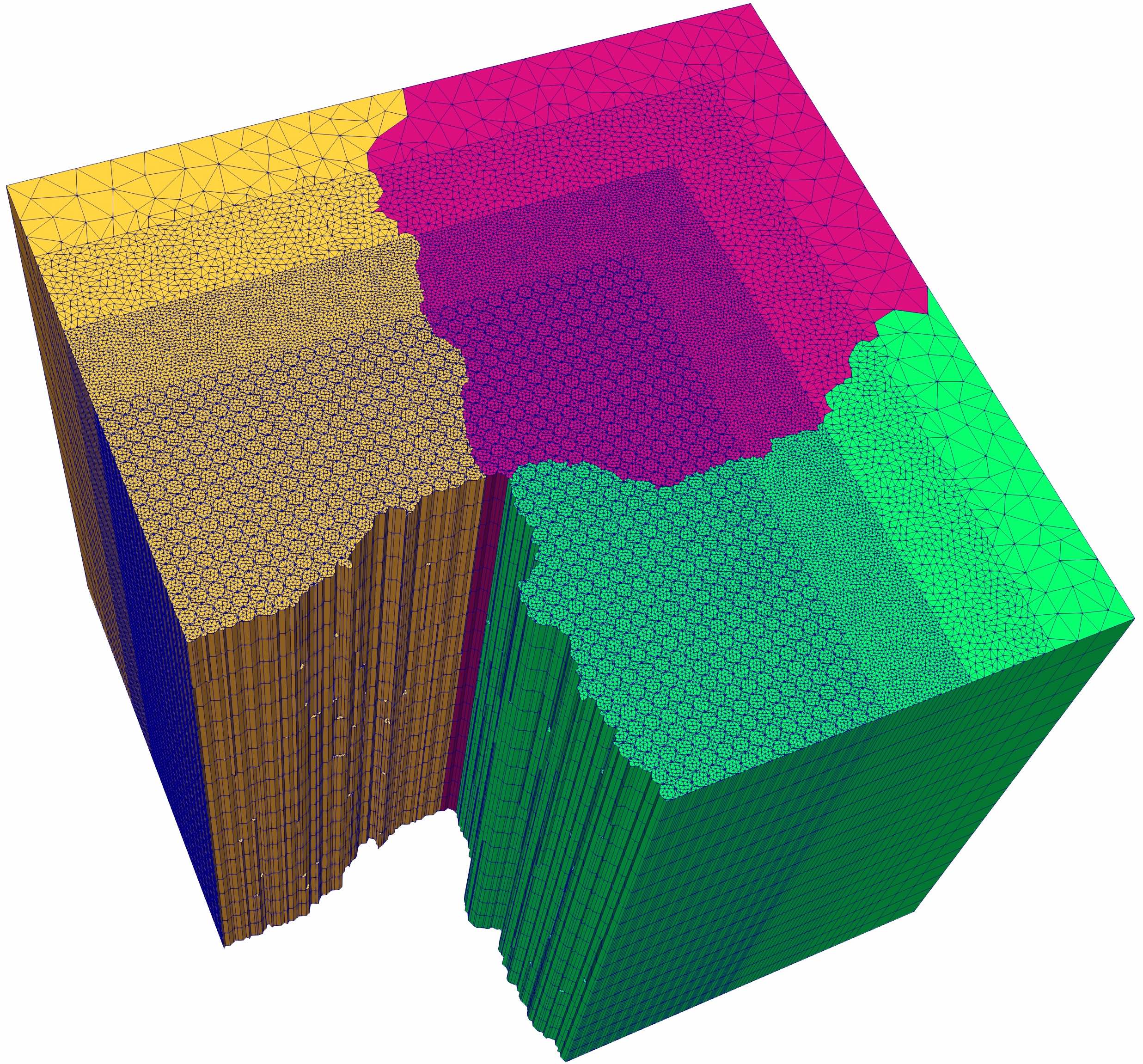} \\
  \includegraphics[width=0.49\linewidth]{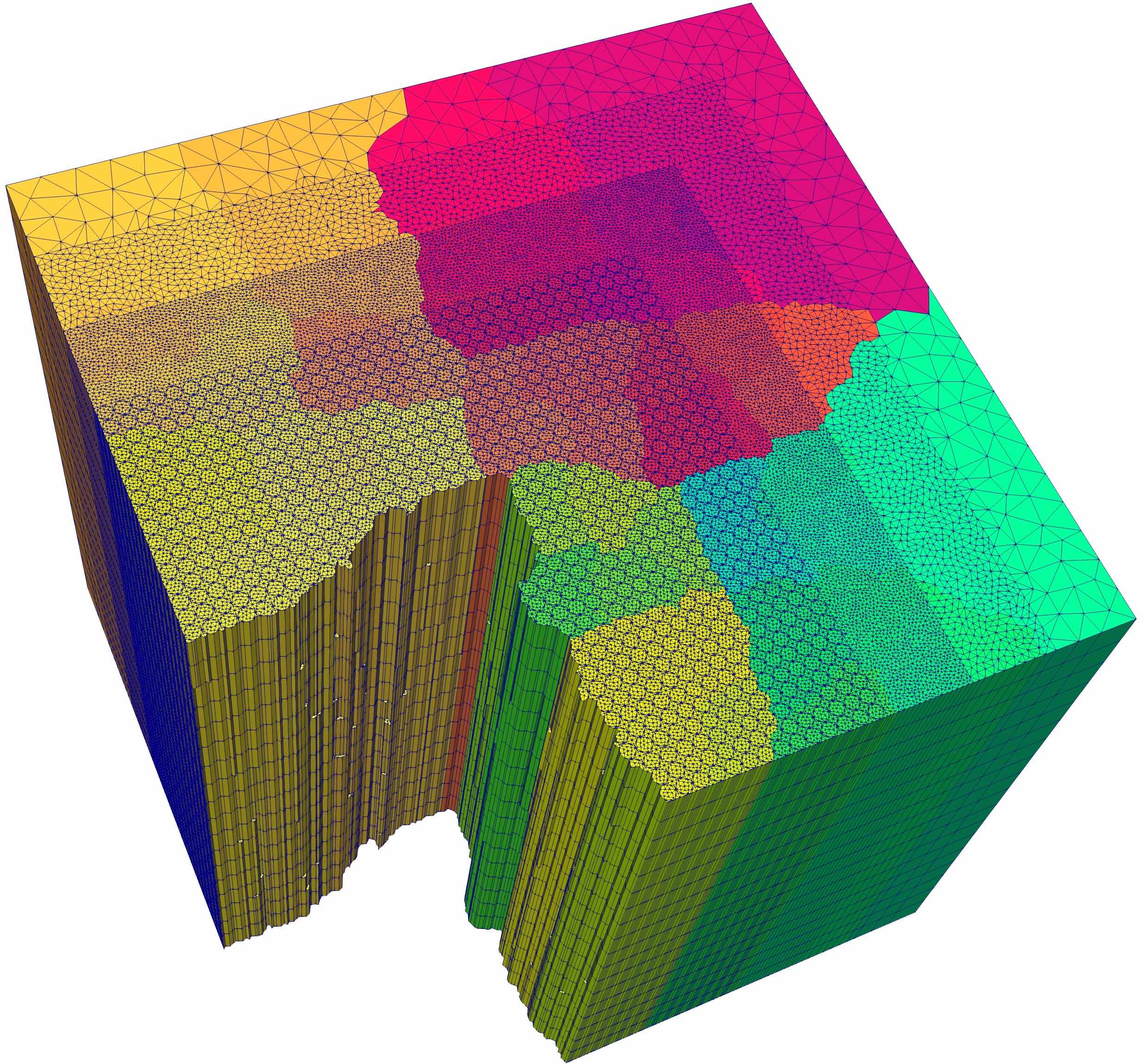} 
  \includegraphics[width=0.49\linewidth]{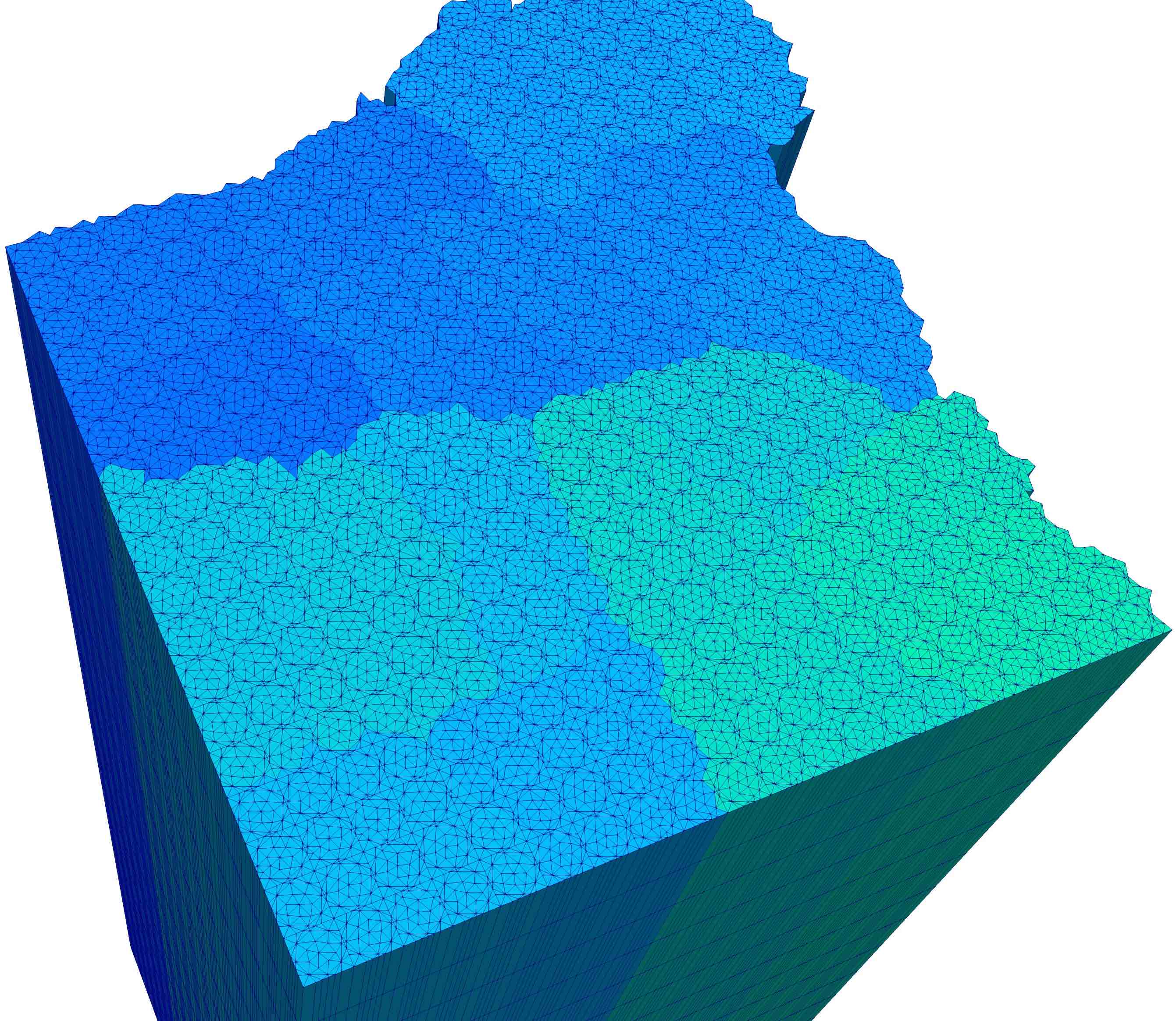} 
 \caption{Demonstration of partitioning a mesh into 80 submeshes using the hierarchical  
 partitioning approach. Assume there 
 are 4 compute nodes available, and each
  compute node has 20 processor cores.
   The mesh is partitioned into 4 submeshes, and each submesh is further partitioned  into 20 small submeshes.
   We have 80 submeshes in total at the end.
   First row: a mesh is partitioned into 4 submeshes; second row: each submesh 
   is divided into 20 small submeshes. For the visualization  purpose,
   the first submesh is taken off in the right picture of the first row. 
   The right picture of the second row is the partition of the first submesh into 20 small submeshes. 
   } \label{fig:hierarch_partition_3D_80}
\end{figure}

In this paper, we employ  a multilevel Schwarz preconditioner 
for the solution of \eqref{alg:newton_preconditioner_process}.
To describe the multilevel preconditioner, 
we partition the fine mesh $\domain_h$, a triangulation of $\domain$,
into $np$ submeshes,$\domain_{h,i}$, $i=1,2,.., np$. Here $np$ is the number of processor cores. 
The mesh is partitioned  using a hierarchical  partitioning  approach \cite{kong2016highly,kong2018general}. 
The hierarchal partitioning takes into consideration that there are multiple processor  cores per each compute 
node on modern supercomputers and these cores share the same chunk of memory.
The hierarchal partitioning generally   consists of two steps 
(even though it has been extended to a 
multilevel version \cite{kong2018general}).  The mesh is first partitioned into $np_1$ 
submeshes ($np_1$ is often the number of compute nodes), 
and then each submesh is further divided  
into $np_2$ small submeshes  ($np_2$ is often  
the number of processor cores per compute node). 
The hierarchical  partitioning  is able to minimize the communication between compute nodes since
only compute-node boundary  cores  need to send messages across the network, and the messages 
within the same compute node are efficiently handled by a modern MPI implementation.  
A detailed discussion on
the hierarchal partitioning can be found in \cite{kong2016highly, kong2018general}.  A hierarchical paritioning
example is shown in Fig.~\ref{fig:hierarch_partition_3D_80}, where a 3D mesh is 
partitioned into 4 submeshes and each 
submesh is further  partitioned into 20 small submeshes. Finally we have 80 small submeshes in total. 

We denote  submatrices  and  subvectors on $\domain_{h,i}$, as $\vM_i$, $\ve_i$ and $\vr_i$, ($i=1,2,\cdots, np$), respectively.
The submarices $\vM_i$ are extended to overlap with their neighbors by  $\delta$ layers, and the resulting  overlapping 
 submatrices are denoted as $\vM_i^{\delta}$. The corresponding subvectors are written as 
 $\ve_i^{\delta}$ and $\vr_i^{\delta}$.  The overlapping-submatrix construction 
 is efficiently implemented  based on 
 the matrix sparsity pattern by defining  
 a restriction operator $\vR_i^{\delta}$ that extracts 
 an overlapping  
 subvector from  the global vector,
 $$
 \vr_i = \vR_i^{\delta} \vr = 
 \left [
 \vI_i^{\delta} ~~ 0
 \right ]
  \left [
 \begin{array}{lll}
 \vr_i^{\delta} \\
 \vr \backslash \vr_i^{\delta}
 \end{array}
  \right ]. 
 $$
Here $\vI_i^{\delta}$ is an identity matrix whose size is the same as that  of $\vr_i^{\delta}$, and
$ \vr \backslash \vr_i^{\delta}$ represents the elements in $\vr$ but not in $\vr_i^{\delta}$.  Similarly, a nonoverlapping version 
of restriction operator is defined as $\vR_i^{0}$ by simply setting $\delta=0$. 
\begin{remark}
A detailed description on the  overlapping submatrix construction
 can be found in \cite{smith2004domain,toselli2006domain}. 
Here, we demonstrate the basic idea using a 1D  problem 
as shown in Fig.~\ref{fig:1Dmesh}. 
\begin{figure}
 \centering
 \includegraphics[width=0.49\linewidth]{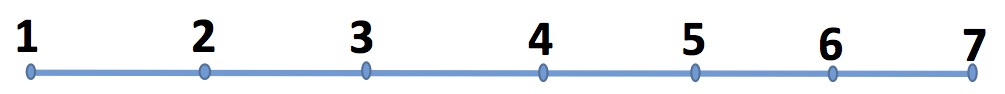}
 \caption{ 1D mesh used for the demonstration 
 of constructing overlapping submatrices.
   } \label{fig:1Dmesh}
\end{figure}
Assume that a one-group diffusion 
equation is discretized  on the 1D mesh using 
the first-oder Lagrange finite element method.  The resulting matrix 
is written as follows
\begin{equation} \label{1Dmatrix}
\vM_{\text{1D}} = 
\left [
\begin{array}{cccccccccc}
m_{11} & m_{12} & &  \\
m_{21} & m_{22} & m_{23} & \\
 & m_{32} & m_{33} & m_{34} \\
  &  & m_{43} & m_{44} &  m_{45} \\
   &  &  & m_{54} & m_{55} &  m_{56} \\
     &  &  & & m_{65} &  m_{66} &   m_{67} \\
     &  &  & & &  m_{76} &   m_{77}   \\
\end{array}
\right ]. 
\end{equation}
It is supposed  that the 1D mesh is partitioned into two parts in such a way that 
vertices 1, 2, 3 and 4 are allocated to processor 1 and vertices 5, 6 and 7 are
owned by processor 2. The submatrices with $\delta=0$ are written as  
\begin{equation} \label{1Dsubmatrix0_1}
(\vM_{\text{1D}})_1^{0} = 
\left [
\begin{array}{cccccccccc}
m_{11} & m_{12} & &  \\
m_{21} & m_{22} & m_{23} & \\
 & m_{32} & m_{33} & m_{34} \\
 &  & m_{43} & m_{44}   \\
 \end{array}
\right ],
\end{equation}
and
\begin{equation} \label{1Dsubmatrix0_2}
(\vM_{\text{1D}})_2^{0} = 
\left [
\begin{array}{cccccccccc}
  m_{55} &  m_{56} \\
m_{65} &  m_{66} &   m_{67} \\
&  m_{76} &   m_{77}   \\
 \end{array}
\right ].
\end{equation}
When $\delta=1$ and $2$, the corresponding overlapping submatrices are written as, respectively,
\begin{equation} \label{1Dsubmatrix1_1}
(\vM_{\text{1D}})_1^{1} = 
\left [
\begin{array}{cccccccccc}
m_{11} & m_{12} & &  \\
m_{21} & m_{22} & m_{23} & \\
 & m_{32} & m_{33} & m_{34} \\
 &  & m_{43} & m_{44}  &  m_{45}\\
  &  & & m_{54} &  m_{55} \\
 \end{array}
\right ],
\end{equation}
\begin{equation} \label{1Dsubmatrix1_2}
(\vM_{\text{1D}})_2^{1} = 
\left [
\begin{array}{cccccccccc}
 m_{44}  &  m_{45}\\
m_{54}&  m_{55} &  m_{56} \\
&m_{65} &  m_{66} &   m_{67} \\
&&  m_{76} &   m_{77}   \\
 \end{array}
\right ],
\end{equation}
\begin{equation} \label{1Dsubmatrix2_1}
(\vM_{\text{1D}})_1^{2} = 
\left [
\begin{array}{cccccccccc}
m_{11} & m_{12} & &  \\
m_{21} & m_{22} & m_{23} & \\
 & m_{32} & m_{33} & m_{34} \\
 &  & m_{43} & m_{44}  &  m_{45}\\
  &  & & m_{54} &  m_{55} &  m_{56}  \\
  &  &  & & m_{65} &  m_{66}  
 \end{array}
\right ],
\end{equation}
and
\begin{equation} \label{1Dsubmatrix2_2}
(\vM_{\text{1D}})_2^{2} = 
\left [
\begin{array}{cccccccccc}
  m_{33} & m_{34} \\
  m_{43}& m_{44}  &  m_{45}\\
&m_{54}&  m_{55} &  m_{56} \\
&&m_{65} &  m_{66} &   m_{67} \\
&&&  m_{76} &   m_{77}   \\
 \end{array}
\right ].
\end{equation}
Note, the 1D mesh is used for the visualization only, 
and the submatrices are constructed algebraically without using a mesh.
\end{remark}

With these notations,
a one-level restricted Schwarz preconditioner 
\cite{kong2018efficient, cai1999restricted} is defined as follows,
\begin{equation} \label{eq:ras}
\vM_{\text{one}}^{-1} =  \sum_{i=1}^{np}(\vR_i^{0})^T (\vM_i^{\delta})^{-1} \vR_i^{\delta},  ~ \vM_i^{\delta} = \vR_i^{\delta} \vM (\vR_i^{\delta})^T,
\end{equation}
where $(\vM_i^{\delta})^{-1}$ represents  the submesh solver 
that is SOR (successive over-relaxation) in this paper.  $\vR_i^{\delta}$  
is not necessarily  implemented 
as a matrix; e.g.,  in PETSc \cite{petsc-user-ref} it is
 implemented as a map from the global vector  to the local  
overlapping subvector. We refer $\vM_{\text{one}}^{-1}$ to as  
a monolithic  preconditioner  since it handles   all unknowns  
 simultaneously, which is fundamentally  different from  
 a traditional  Gauss-Seidel approach.   
 Different  variants of domain composition methods can be found 
in other  literatures \cite{smith2004domain,toselli2006domain}. 

$\vM_{\text{one}}^{-1}$ works for many problems, 
but  coarse spaces are often required 
to enhance the preconditioner  
when the number of processor cores is large  
or when  the system is ill-conditioned. 
In general, there are two approaches 
to generate coarse spaces. The first approach  
is to coarsen the fine mesh to generate coarse finite element 
meshes based on which  the coarse spaces are constructed. 
This geometric method has been
successfully  applied  to many applications; e.g., 
elasticity  problems \cite{kong2016highly}, 
fluid-structure interactions \cite{kong2017scalable, kong2018scalability,kong2016parallel, kong2019simulation}.  
The other approach is to avoid coarsening 
the fine mesh, instead, coarse spaces are generated using the matrix  only.  
The first approach is referred to as ``geometric multilevel 
method", and the second approach is referred  to as ``algebraic multilevel method".  
We consider the algebraic  approach here
since it has an out-of-box  nature. 
Note that no coarse finite element meshes are constructed  in the algebraic
multilevel method, and the required coarse spaces are built algebraically  with a matrix coarsening  algorithm 
to be discussed shortly.
However, algebraically constructing coarse spaces is expensive 
and sometimes  unscalable.  In order to make the overall algorithm scalable, 
we explore the matrix structures  to reduce the cost of coarse space 
 construction.  
 
 For saving memory, the preconditioning matrix ignores the energy coupling in the 
 diffusion system, and drops the energy and angular coupling in the transport system.
 If the matrix 
 was  ordered component-by-component, 
 we  would have a block diagonal matrix,
 \begin{equation} \label{matrix-structure}
 \vM = 
 \left [
 \begin{array}{lllll}
 M_{1} &  & & \\
 &  M_{2} &  & \\
 &&  \ddots& \\
 &  &  & M_{n_{\text{comp}}}  \\ 
 \end{array}
 \right ],
 \end{equation}
 where $n_{\text{comp}}$ is the number of components.
In the low order diffusion system, each component corresponds to an energy group,
 that is, $n_{\text{comp}} =G$.
 In the transport system, each component corresponds to an angular direction of a given energy 
 group, that is, $n_{\text{comp}} = G \times N_{\text{d}}$.
 In either the diffusion system or the transport system, $M_{j}, j=1, 2, \cdots, n_{\text{comp}}$,
 is a spatial matrix.
 
Component matrices (spatial matrices), $M_{j}, j=1, 2, \cdots, n_{\text{comp}}$, 
  correspond  to the discretization  of  partial differential   equations  
 on $\domain_h$.
 The structures of spatial matrices are similar to 
 each other since they  
 correspond  to the same partial differential   
 equation operators, and the numeric 
 values of spatial matrices  
 are different from each other
 because of different materials. 
 Based on these observations,
 we coarsen  a single-component matrix, e.g., $M_{1}$, instead 
 of the full matrix, $\vM$, to generate a sequence of  subinterpolations 
 that are reused for all other components. The subinterpolations
 generated based on a single-component matrix are referred to as
 ``Single-component Grid based interpolations" (SG interpolations). More precisely,  
 we use GAMG \cite{petsc-user-ref} or BoomerAMG \cite{yang2002boomeramg} to 
 coarsen $M_{1}$ to build  a $L$-level hierarchy   
 consisting  of $(L-1)$ subinterpolations, $P_{l+1}^l$, $l=1,2,.., L-1$.  
 The coarsening process based on a single-component matrix 
 to generate a sequence of SG interpolations  is referred to 
 as ``subspace-based coarsening".
 To define the full interpolations, 
 we introduce a restriction operator 
 that extracts a subvector $r_j$ for component $j$  from the full vector $\vr$,
 $$
 r_j =  R_j  \vr = 
 \left [
 I_j ~~ 0
 \right ]
  \left [
 \begin{array}{lll}
 r_j \\
 \vr \backslash r_j
 \end{array}
  \right ], j=1,2, \cdots, n_{\text{comp}}. 
 $$
 Here $I_j$ is an identity matrix, and $\vr \backslash r_j$ denotes 
 all the elements in $\vr$ but not in $r_j$. $R_j$ can be defined 
 on all levels, and it is rewritten as $R_j^{l}$ for the $j$th component on the $l$th level.
 The full interpolations are then written as
 \begin{equation}\label{asm:fullinterpolation}
 \vP_{l+1}^l =  \sum_{j=1}^{n_{\text{comp}}} (R_j^{l})^T P^l_{l+1} R_j^{l+1}, l =1, 2, ..., L-1.
\end{equation}
 The coarse operators are defined using the Galerkin method with  $ \vP_{l+1}^l$, that is,
 \begin{equation}\label{asm:operators}
 \vM^{l+1} =  (\vP_{l+1}^l )^T \vM^{l} \vP_{l+1}^l.
\end{equation}
Here $\vM^1 = \vM$ is the finest level operator.  The coarse spaces defined in~\eqref{asm:operators} are
referred to as ``Single-component Grid based coarse spaces" (SG coarse spaces) since they are 
constructed based on  the SG interpolations.
 
Similarly, the one-level method can be defined on each level, 
and these one-level preconditioners 
are denoted as $(\vM_{\text{one}}^{l})^{-1}$ ($l=1,2,\cdots, L$).  
 The multilevel method constructed 
 with the subspace-based coarsening algorithm 
 is referred to as 
 ``Single-component Grid based Multilevel Additive Schwarz Method" (SGMASM), 
 shown in Alg.~\ref{SGMASM}. The method using the unmodified standard
 coarse spaces is denoted as ``Multilevel Additive Schwarz Method" (MASM), 
 shown in Alg~\ref{MASM}.   For saving memory, 
 the interpolations in~\eqref{asm:fullinterpolation}
are implemented without storing duplicated 
values in a specific matrix format in 
PETSc \cite{petsc-user-ref}, and
the sparse matrix triple products in~\eqref{asm:operators} 
are carried out using an all-at-once approach \cite{kong2019parallel}. 
Note, in Alg.~\ref{SGMASM}, the SG coarse spaces are generated 
based on the first component matrix. In practice,
any component matrix can be chosen to generate a sequence of SG coarse spaces, 
and it is an user-specified parameter in  software.
\begin{algorithm}
\caption{Application of multilevel preconditioner}\label{MASM-application}
\begin{algorithmic}[1]
\Procedure {PCApply} {$\vM^{l}, \ve^{l}, \vr^{l}$}
\If {$l=L$} 
 \State Solve  $\vM^{L} \ve^{L} = \vr^{L}$ with a  direct solver redundantly on each compute node or each processor core
\Else 
 \State Pre-solve $\vM^{l} \ve^{l} = \vr^{l}$ for $\ve^{l}$ using a $(\vM^{l}_{\text{one}})^{-1}$ preconditioned iterative method
 \State  Set $\bar{\vr}^{l} = \vr^{l} - \vM^{l} \ve^{l}$
 \State  Apply  restriction: $ \bar{\vr}^{l+1} = (\vP_{l+1}^l)^T \bar{\vr}^{l}$
 \State  Call PCApply($\vM^{l+1}, \vz^{l+1},  \bar{\vr}^{l+1}$)
 \State  Apply  interpolation: $\vz^{l} = \vP_{l+1}^l \vz^{l+1}$
 \State  Apply  correction: $\ve^{l} = \ve^{l} + \vz^{l}$
 \State Post-solve  $\vM^{l} \ve^{l} = \vr^{l}$ for $\ve^{l}$ using a  $(\vM^{l}_{\text{one}})^{-1}$ preconditioned iterative method
\EndIf
\EndProcedure
 \end{algorithmic}
\end{algorithm}
 \begin{algorithm}
\caption{Multilevel  Schwarz preconditioner constructed  with the subspace-based coarsening algorithm (SGMASM)}\label{SGMASM}
\begin{algorithmic}[1]
\Procedure {PCSetup} {$\vM$} 
\State  Extract a component matrix  $M_{1}$ from $\vM$ in parallel
\State  Coarsen $M_{1}$ and  generate ($L-1$) subinterpolations $P_{l+1}^l$
\For {$l=1,2,...,L-1$}
\State  Construct  the $l$th full interpolation:   $\displaystyle  \vP_{l+1}^l = \sum_{j=1}^{n_{\text{comp}}}  (R_{j}^{l})^T P_{l+1}^l  R_{j}^{l+1}$
\EndFor
\For {$l=1,2,..., L-1$}
\State  Build  the $(l+1)$th coarse matrix: ${\vM}^{l+1} =   (\vP^{l}_{l+1})^T \vM^{l}  \vP_{l+1}^l$
\EndFor
\State Return $\{ {\vM}^{l+1} \}$ and $\{ \vP_{l+1}^l \}$, $l=1,2,..,L-1$
\EndProcedure
\State Call PCApply($\vM^{1}, \ve^{1}, \vr^{1}$)  \Comment{ Alg.~\ref{MASM-application}}
\State Return $\ve^1$
 \end{algorithmic}
\end{algorithm}
 \begin{algorithm}
\caption{Traditional multilevel  Schwarz preconditioner  (MASM)}\label{MASM}
\begin{algorithmic}[1]
\Procedure {PCSetup} {$\vM$} 
\State  Coarsen $\vM$ and  generate ($L-1$) interpolations $\vP_{l+1}^l$ and ($L-1$) coarse matrices $\vM^{l+1}$
\State Return $\{ {\vM}^{l+1} \}$ and $\{ {\vP}_{l+1}^l \}$, $l=1,2,..., L-1$
\EndProcedure
\State Call PCApply($\vM^{1}, \ve^{1}, \vr^{1}$)  \Comment{ Alg.~\ref{MASM-application}}
\State Return $\ve^1$
 \end{algorithmic}
\end{algorithm}

\section{Numerical results}
\begin{figure}
 \centering
  \includegraphics[width=0.51\linewidth]{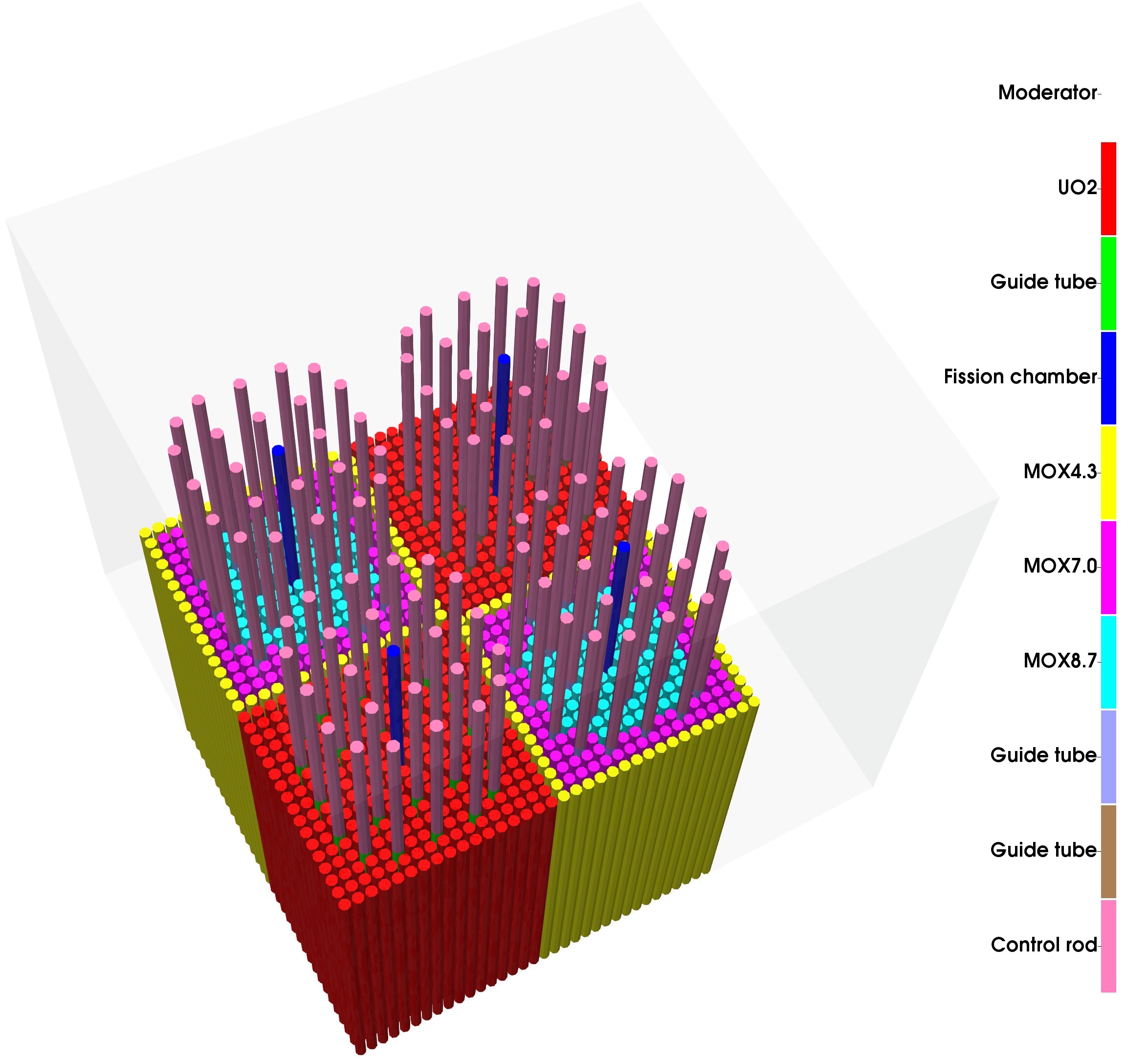} 
  \includegraphics[width=0.46\linewidth]{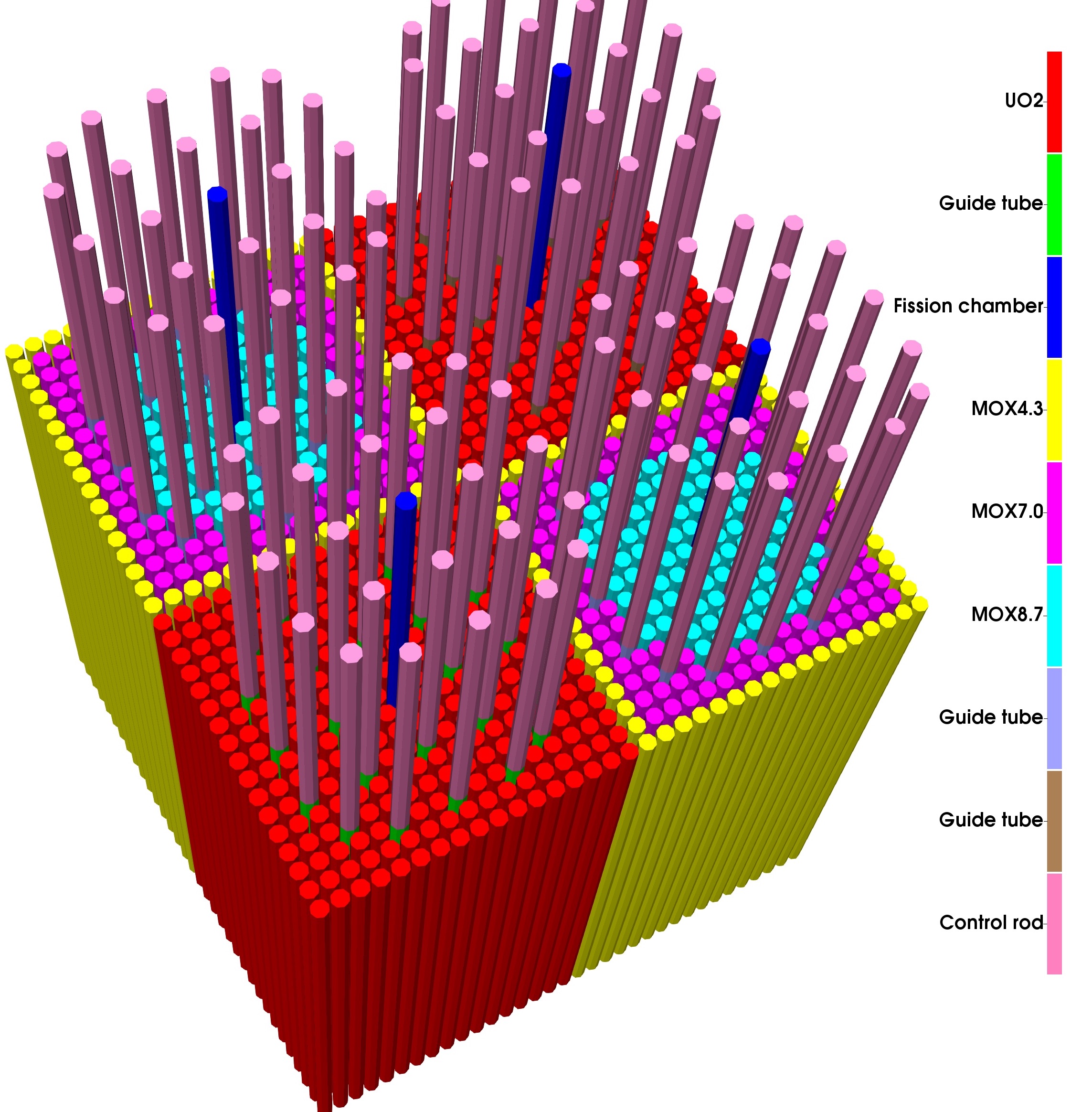} 
 \caption{C5G7 3D benchmark configuration. 
 Left: four assemblies  at the front bottom  
 of  a $64.26 \times 64.26 \times 64.26$ cm cube; 
 right: four assemblies and control rods. 
 Different pin  colors  correspond to different  materials, 
 e.g., UO2, guide tube, fission chamber, MOX 4.3\%,
MOX 7.0\%, MOX 8.7\%, and control rod. } \label{fig:configuration_3D_fuel}
\end{figure}

\begin{figure}
 \centering
  \includegraphics[width=0.49\linewidth]{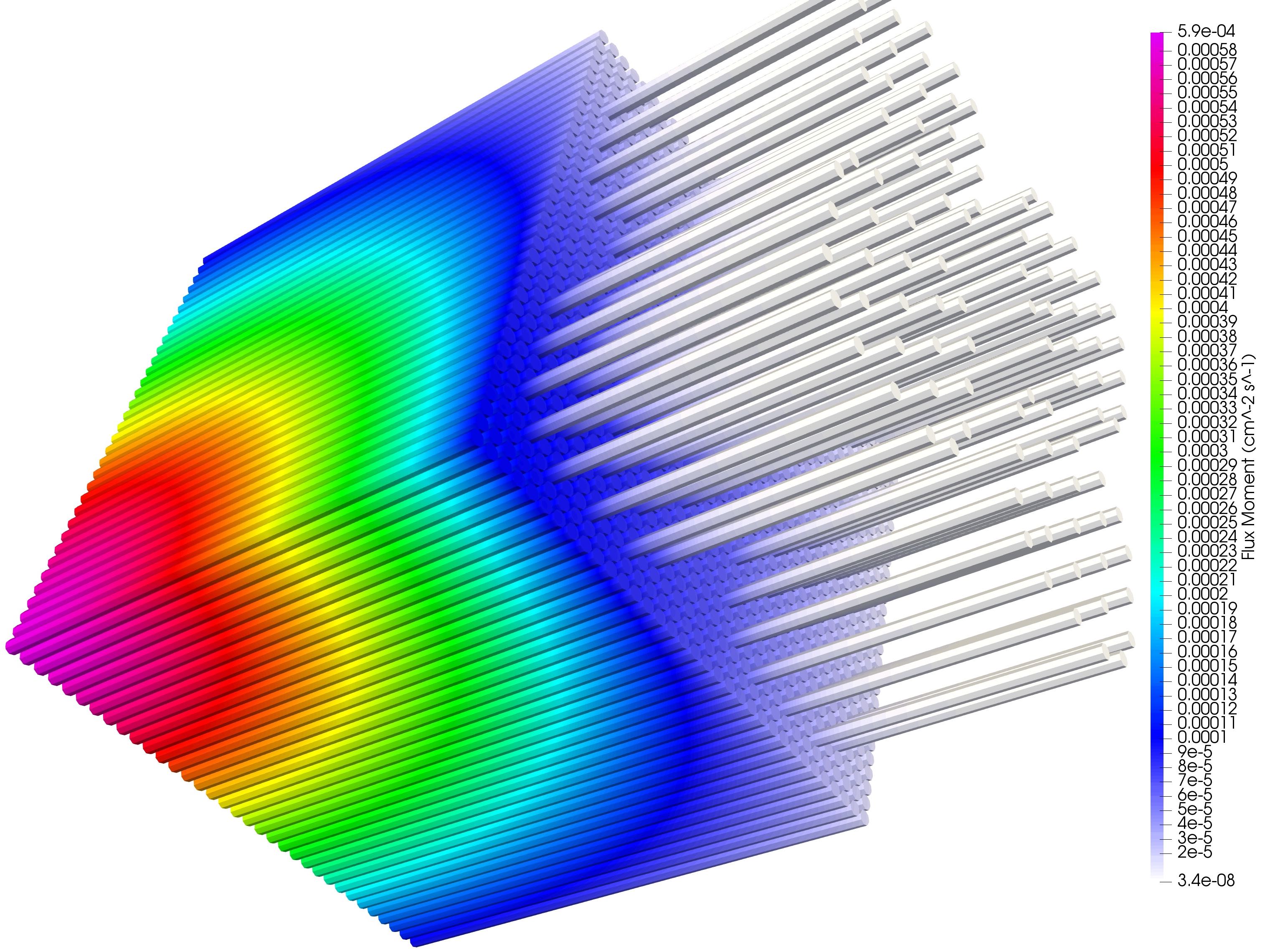} 
  \includegraphics[width=0.49\linewidth]{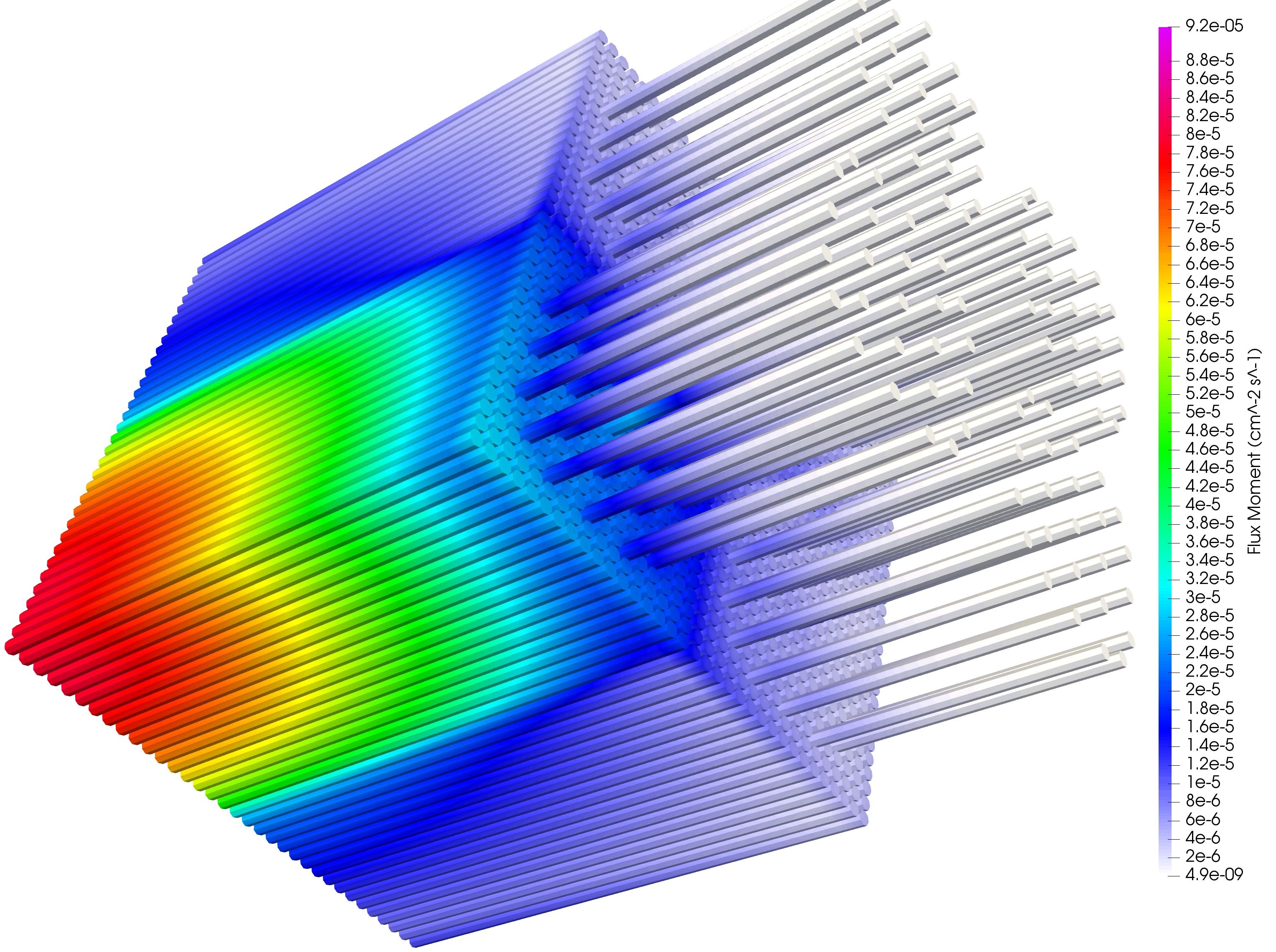} \\
 \caption{Flux moments. Left: the first group flux moment; right: the sixth group flux moment.} \label{fig:solution}
\end{figure}

In this section, we  report the performance of the proposed algorithm 
 in terms of the compute time 
and the memory usage.  The C5G7 MOX 3D benchmark is employed  
here as a test case, and its configuration 
is shown in Fig.~\ref{fig:configuration_3D_fuel}, 
where different  pin  colors  correspond 
to different  materials; e.g., UO2, guide tube, fission chamber, MOX 4.3\%,
MOX 7.0\%, MOX 8.7\%, and control rod.   The overall dimensions of the domain   
  are $64.26 \times  64.26 \times 64.26$ cm, 
where  each  fuel  assembly is $21.42 \times 21.42 \times 42.84$ cm.  A fuel assembly 
 consists  of a $17 \times 17$ lattice of  pin cells, 
where    the side length of each pin cell is 1$.26$ cm and the radius  of the fuel pins and guide 
tubes is  $0.54$ cm. 
The reflected boundary conditions are applied to the front, the left and the bottom boundaries, and 
 the vacuum boundary conditions are applied to the back, the right and the top boundaries.
The seven-group 
set of cross sections are often chosen to test the algorithm performance  since the corresponding 
problem is difficult to solve.   More details of the benchmark configuration 
can be found in literatures \cite{kosaka20063, lewis2001benchmark}.   
The computed eigenvalue is 1.141932, and the eigen flux moments  for the second 
and the sixth groups are shown in Fig.~\ref{fig:solution}.  

The preconditioner is implemented 
in PETSc \cite{petsc-user-ref} as part of this work,
and the subspace-based coarsening algorithm is implemented 
based on BoomerAMG \cite{yang2002boomeramg}.
The angular and spatial discretizations  are implemented 
in Rattlesnake \cite{wang2018rattlesnake} that is 
on  top of MOOSE \cite{peterson2018overview,permann2019moose} and libMesh \cite{libMeshPaper}.  
The numerical experiments  
are carried out on a supercomputer at INL (Idaho National Laboratory), 
where each compute node has two 20-core processors with 2.4 GHz  
and the compute nodes  are connected by an OmniPath  network.

For simplifying the discussion, we define some notations 
that will be used in the rest of paper.  
``Mem" represents 
the estimated memory usage per processor core in Megabyte (M),  
``Its$_{\text{Newton}}$" 
is the averaged number of Newton iterations per   
Picard iteration for the diffusion system,  ``Its$_{\text{linear}}$"  
is the averaged number of GMRES iterations 
per Newton step in the diffusion system, 
``Its$_{\text{sweep}}$" denotes the averaged  
number of GMRES iterations per Picard iteration  
in the transport system,  
``Time$_{\text{low}}$" represents  the compute time 
spent on the diffusion solver,  ``Time$_\text{high}$" 
is the compute time on the transport  solver, ``Time$_T$" is 
the total compute time for the overall simulation, 
and ``EFF" is the parallel efficiency 
with respect to the number of processor cores. 
``Time$_\text{ksp}$"   is the total compute time on the 
linear solver,  ``Time$_{\text{PCA}}$" represents 
the total compute  time on the preconditioner application,  
``Time$_{\text{PCS}}$"  is 
the total compute  time on the preconditioner setup,  
``Time$_{\text{MF}}$" denotes the compute time 
on the matrix-free matrix-vector operations, 
``Time$_{\text{Func}}$" is the time spent on the function evaluations, 
``Time$_\text{Jac}$" represents  the compute time  on the Jacobian evaluations,
and ``Time$_\text{LS}$"  is the time consumed  by the line search routine. 
All compute times are reported in second (s).
A relative tolerance of $10^{-5}$ is chosen for the linear solver of   the transport system.
In the diffusion system, Newton is stopped when a relative tolerance of $10^{-3}$ is met, and 
the relative tolerance  of the linear solver is $10^{-2}$.
Below we start to discuss   the algorithm robustness with  
respect to different parameters; e.g., subdomain overlapping size, strength matrix threshold and 
the number of aggressive coarsening levels. 

\subsection{Subdomain overlapping size}
A subdomain overlapping size plays an 
important  role on domain decomposition methods.  
A larger overlapping size often 
leads to a better convergence in terms of the number 
of GMRES iterations, meanwhile the resulting preconditioner  also
consumes more memory and  involves more communication. 
We need to take the iteration count,
the memory usage and the communication 
into consideration to choose an optimal overlapping size, 
and the optimal overlapping size   is often problem dependent.   

In this test,
we report the parallel  performance of the proposed algorithm  
with respect to different overlapping sizes. 
A Gauss Chebyshev angular quadrature scheme with 32  directions is employed.  
A mesh with 832,371 nodes and  1,567,944 elements
is used for both the diffusion system and the transport system.
The diffusion system has 5,826,597 unknowns, and 
the transport system has 186,451,104 unknowns.  The numerical 
results are summarized in Table~\ref{tab:overlap}.   
The neutron transport problem is accelerated by the low order diffusion 
system. For solving the linear system of equations 
in the transport system, the preconditioning matrix 
is coarsened with the subspace-based coarsening 
algorithm to generate 7 subinterpolations, 
and a seven-level preconditioner 
is constructed using these subinterpolations. 
On the finest level a Schwarz preconditioner 
with different overlapping sizes is employed, and no overlapping is
used on all other coarse levels except the coarsest level on which  
a direct solver is carried out  redundantly.   The low order 
diffusion system is computed  with Newton-Krylov-SGMASM 
with an initial guess obtained  from two
inverse power iterations.
The linear system of transport  equations is also solved 
  with SGMASM preconditioned GMRES 
  with a seven-level 
hierarchy. 
\begin{table}
\scriptsize
\centering
\caption{Parallel performance with respect to different subdomain overlapping sizes.  The transport system with 
 186,451,104 unknowns is accelerated  by the low order diffusion  system with 5,826,597 unknowns.  \label{tab:overlap}}
\begin{tabular}{c c c c c c  c c c c}
\toprule
$np$ & $\delta$ &Mem (M) & Its$_{\text{newton}}$ & Its$_{\text{linear}}$ &Its$_{\text{sweep}}$
&Time$_{\text{low}} (s)$&Time$_\text{high}$ (s)& Time$_T$ (s)& EFF  \\
\midrule
160 &0 &  1587 &  2 & 32&30&661& 626 &1287 &100\%  \\
160 & 1  &  2302 &2 &28&28 & 636& 653& 1289&100\% \\
160 & 2  &  2680 &2 &26&25& 595&  678& 1273&100\% \\
\midrule
320 &0 &  810 &  2 & 33&30&348&  315&  663& 97\%  \\
320& 1 &  1194 &2 &28&28 &333&  336& 669& 96\% \\
320 & 2  &  1422 &2 &26&25& 314& 370& 684& 93\% \\
\midrule
640 &0  &  427 &  2 & 34&30&187 &165&352&  91\%  \\
640 & 1  &  656 &2 &29&29 & 185&187& 372& 87\% \\
640 & 2  &  798 &2 &26&26& 171&220& 391 &81\% \\
\midrule
1,280& 0 &  210 &  2 & 35&31&111&89&200& 80\%  \\
1,280 & 1  &  341 &2 &29&28 &106&108&214& 75\% \\
1,280 & 2  &  526 &2 &26&26&103& 136&239&67\% \\
\bottomrule
\end{tabular}
\end{table}

In Table~\ref{tab:overlap}, it is easily observed  
that more memory is used as we increase 
the subdomain  overlapping size, especially  
when $\delta$ increases  from 0 to 1.  ``$\delta=1$" 
uses much more memory than ``$\delta=0$" because 
the local submatrix   with ``$\delta=0$" is  
implemented in place while the implementation of  ``$\delta=1$"
has to involve an extra submatrix copy 
for storing overlapping and local  elements. 
For the 160-core case,  the overall algorithm with $\delta=0$ 
consumes $1587$ M memory (per processor core),
and the memory usage is increased to $2302$ M 
by $715$ M when we increase $\delta$ to 1. 
The memory usage continues being increased to $2680$ M 
by $378$ M when we continue increasing $\delta$
to 2.  The number of Newton iterations per Picard iteration 
stays as a constant for all overlapping sizes and all core counts.
The number of GMRES iterations per Newton step   
in the diffusion system  decreases slightly 
when a larger overlapping size is employed so that the compute time on the diffusion system 
 is reduced accordingly. 
GMRES iteration for the transport system  is, similarly, 
decreased when we change $\delta$ from 0,  to 1 and 2, 
but the compute time on the transport 
system  does not decrease, instead,  it is increased from $626$ s, to $653$ s and $678$ s 
when 160 processor 
cores are used because the per-iteration cost with 
a larger  $\delta$ is
higher than   that obtained with a smaller   $\delta$. 
That is,  the reduction in GMRES iteration can not compensate the per-iteration-cost increase.  
The total 
compute times with different overlapping sizes at 160 processor cores are close to each other
because the time reduction  in the diffusion system is cancelled 
by the time increase in the transport system.
 
\begin{figure}
 \centering
 \includegraphics[width=0.49\linewidth]{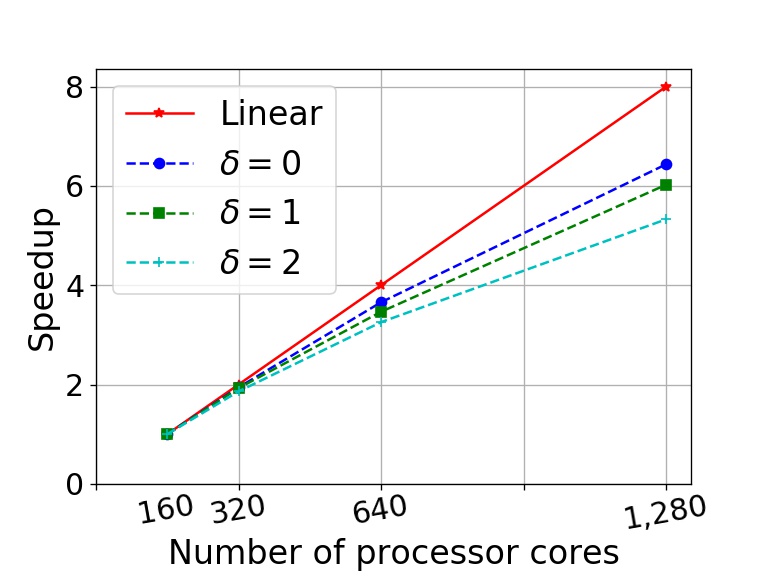} 
 \includegraphics[width=0.49\linewidth]{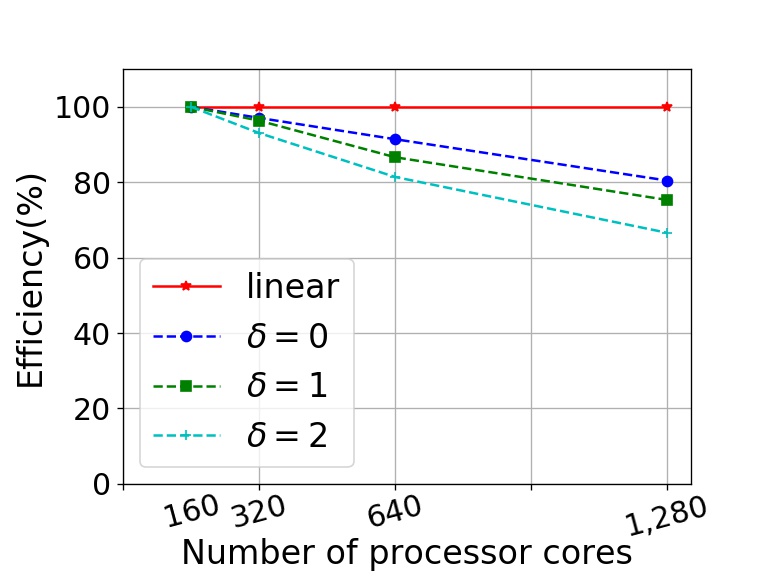} 
 \caption{Speedup and efficiency for different overlapping sizes on 160, 320, 640 
 and 1,280 processor cores. Left: speedup, right: parallel efficiency.\label{fig:overlap-speedups}}
\end{figure}

We now try to understand how the  overall algorithm performance correlates to  the 
overlapping size when the number of processor cores is increased.
The proposed algorithm  with $\delta=0$ is scalable in memory  in the sense
that the memory usage is halved when we double the number of processor cores.
More precisely, the memory usage with $\delta=0$ is almost halved to 
$810$ M from $1587$ M when the core count is doubled  from 
160 to 320.  It continues being reduced to $427$ M and $210$ M, respectively, 
when the number of processor cores is  increased to $640$ and $1,280$.
For both $\delta=1$ and $2$,  we observed the similar behaviors, that is,
the overall algorithm equipped with  different overlapping sizes is scalable in  memory.
The memory usages with $\delta=1 \text{~and~} 2$  are much higher than 
that consumed with $\delta=0$ for all core counts. At 1,280 processor cores, 
the memory consumed with $\delta=1$ is almost twice as much as that using $\delta=0$.
The memory usage  using $\delta=2$ is  $2.5$ times as much as that consumed with  $\delta=0$. 
The overall algorithm is mathematically scalable since the 
averaged numbers  of GMRES iterations  
in both the transport  and the diffusion systems stay close to constants  for all 
core counts.  The compute times in the diffusion system with $\delta=0$
are more than other cases for all core  counts, while the algorithm performs
better with a smaller overlapping size  in the transport system, e.g., $\delta=0$
is better than the nonzero overlapping sizes.   
We conclude that  $\delta=0$ is an optimal choice for the transport system,
and  $\delta=2$ is a better choice for the diffusion system. 
The overall algorithm is highly scalable for different overlapping sizes with up to
1,280 processor cores. Parallel efficiencies of $80\%$, $75\%$ and $67\%$ are obtained when 
using $\delta=0, 1  \text{~and~} 2$ at 1,280 processor cores. The speedups  
and parallel efficiencies 
are also  shown  in Fig.~\ref{fig:overlap-speedups}, where the algorithm using $\delta=0$
is obviously better than all other cases.  

To further understand the preconditioner performance, the compute times on the individual  
 components of  the diffusion system are reported 
in Table~\ref{tab:overlap-low-order}, where the compute times on 
the preconditioner setup and application increase for all core counts 
\begin{table}
\scriptsize
\centering
\caption{Compute times on the individual  components  of  the diffusion system using  different overlapping sizes on 160,
320, 640 and 1,280 processor  cores.  \label{tab:overlap-low-order}}
\begin{tabular}{c c c c c c  c c c c}
\toprule
$np$ & $\delta$ &Time$_{\text{ksp}}$ (s) & Time$_{\text{PCA}}$ (s)
& Time$_{\text{PCS}}$ (s) &Time$_{\text{MF}}$ (s)&Time$_{\text{Func}}$ (s)
&Time$_\text{Jac}$ (s)& Time$_\text{LS}$ (s) \\
\midrule
160 &0 &  584 &  38 & 15&547&28& 31 &37   \\
160 & 1  &  559 &45 &17&514 & 28& 31& 37 \\
160 & 2  &  518 &48 &18&470& 28&  31& 36 \\
\midrule
320 &0 &  307 &  20 & 9&286&15&  16&  19  \\
320& 1 &  293 &25 &10&267 &14&  16& 19 \\
320 & 2  &  274 &28 &11&245& 14& 16& 19 \\
\midrule
640 &0  &  167 &  13 & 6&151&8 &8&10 \\
640 & 1  &  164 &16 &6&145 & 7&8& 10 \\
640 & 2  &  150 &19 &7&130& 7&8& 10 \\
\midrule
1,280& 0 &  99 &  10 & 6&86&4&4&6  \\
1,280 & 1  &  94 &13 &6&75 &4&5&5 \\
1,280 & 2  &  92 &14 &5&75&4& 4&5 \\
\bottomrule
\end{tabular}
\end{table}
when the  overlapping size is increased. The linear solver time is 
reduced when we increase the overlapping size because
the corresponding GMRES iteration is decreased. For example, at 160 processor cores, 
the linear solver time is decreased 
from $584$ s to $559$ s and $518$ s when $\delta$ is 
increased from 0  to 1 and 2.  The compute time on the matrix-free
matrix-vector operations is also decreased when we increase the overlapping size, 
which is observed for all  
core counts.  For instance, at 640 processor cores, the compute time on 
the matrix-free matrix-vector operations
is reduced from $151$ s  to $145$ s when $\delta$ is increased from 0 to 1, 
and it continues being decreased  to $130$ s
by  $15$ s with  using $\delta=2$.   The compute times on the Jacobian evaluations, 
the function evaluations and the line search 
are almost the same for different overlapping sizes since the number of Newton 
iterations stay the same regardless of the overlapping size.
In the diffusion system, the matrix-free matrix-vector operations dominate the 
whole calculations so that a  reduction on the matrix-free
operations due to a large overlapping size  leads to a more efficient computation. 

In the transport system,  the scattering and the fission terms are computed 
using the scalar fluxes 
and the eigenvalue from the diffusion system, and then 
only a linear system of equations needs  to be solved at each Picard iteration.  
The compute times for the individual components of the transport
system are shown in Table~\ref{tab:overlap-high-order}.
\begin{table}
\scriptsize
\centering
\caption{Preconditioner performance in the transport system with 
different overlapping sizes on 160, 320, 640 and 1,280 
processor cores.  \label{tab:overlap-high-order}}
\begin{tabular}{c c c c c c  c c c c}
\toprule
$np$ & $\delta$ & Time$_{\text{PCS}}$ (s) & Time$_{\text{PCA}}$ (s)  & Time$_\text{high}$ (s)  \\
\midrule
160 &0 &  18 &  549 & 626   \\
160 & 1  &  21 &582 &653 \\
160 & 2  &  24 &611 &678 \\
\midrule
320 &0 &  10 &  272 & 315  \\
320& 1 &  12 &298 &336 \\
320 & 2  &  14 &335 &336 \\
\midrule
640 &0  &  7 &  140 & 352 \\
640 & 1  &  7 &166 &372 \\
640 & 2  &  9 &199 &391 \\
\midrule
1,280& 0 &  4 &  73 & 200  \\
1,280 & 1  &  5 &95 &214 \\
1,280 & 2  &  6 &123 &239 \\
\bottomrule
\end{tabular}
\end{table}
We observed that the preconditioner setup time is negligibly  
increased when using different $\delta$, and it takes only a 
small portion of the total compute time because the preconditioner 
is fixed in the transport system for the entire simulation so that 
only one preconditioner setup is involved at the beginning  of the simulation.  
Most of the compute time in the transport system 
is spent on applying the preconditoner to the linear system of equations.  When we increase the overlapping size,
the preconditioner application time is increased for all core counts. For example, at 160 processor cores,
the preconditioner application time  is increased from $549$ s to $582$ s by $33$ s when $\delta$ is increased 
from 0 to 1. It continues  being increased  to $611$ s when  $\delta=2$.  At 1,280 processor cores,
the preconditioner application  times are $73$ s, $95$ s and $123$  s for $\delta=0, 1  \text{~and~} 2$, respectively.  The preconditioner 
application time at $\delta=2$ is almost twice as much as  that at $\delta=1$.  The overall compute time  with $\delta=0$, in 
the transport system, is better than the other two options.

\subsection{Comparison with an unaccelerated transport solver}
The multigroup neutron  transport equations  can be solved directly 
without a nonlinear diffusion acceleration
method (we have done this in our previous work 
\cite{kong2019scalable, kong2019highly}), where 
 a  generalized  eigenvalue problem  instead of 
 a linear system of equations is involved.  
 The goal of the nonlinear diffusion acceleration method   is to 
 decouple the fission and the scattering terms from other terms 
 in the transport equations so that 
 only a linear system of equations needs to be computed, 
 which in turn reduces  the computational  cost.
 In this test, we do a performance comparison between  
 the unaccelerated transport solver and the diffusion 
 accelerated  transport method.  In the diffusion 
 accelerated transport method, ``$\delta=0$" is 
employed for the transport system and ``$\delta=2$" 
is adopted for the diffusion system. 
In the unaccelerated transport solver, 
an inexact Newton-Krylov method together with 
SGMASM is used to directly solve the generalized eigenvalue problem. 
More details on the unaccelerated transport  solver 
can be found in our previous work \cite{kong2019highly,kong2018fully}. 

The performance comparison between the unaccelerated 
transport solver and the diffusion accelerated  transport method 
is reported in Table~\ref{tab:direct-vs-nda}. 
The test is carried out using 1,280 processor cores with different 
numbers  of angular directions.  ``NDA" represents 
the nonlinear diffusion accelerated transport solver, and
``Unaccelerated" denotes the unaccelerated transport solver. 
\begin{table}
\scriptsize
\centering
\caption{Performance comparison between the unaccelerated 
transport solver and the diffusion accelerated 
 transport method  using 1,280 processor cores.  \label{tab:direct-vs-nda}}
\begin{tabular}{c c c c c c  c c c c}
\toprule
$N_\text{d}$ & algrotihm  &Mem (M) & Time$_{\text{MF}}$ (s) & Time$_{\text{PCS}}$ (s) 
&Time$_{\text{PCA}}$ (s)&Time$_{\text{ksp}}$ (s)& Time$_T$ (s)  \\
\midrule
32 & NDA &  212 &  77 & 10&88&179& 194  \\
32 & Unaccelerated  &  256 &246 &9&28 & 283& 313 \\
\midrule
64 &NDA &  433 &  71 & 14&165&261&  281  \\
64&Unaccelerated &  496 &567&16 &56&639 &  705 \\
\midrule
128 &NDA  &  868 &  72 & 20&353&476 &502 \\
128 & Unaccelerated  &  985 &1241 &29&113 & 1382&1527 \\
\midrule
192& NDA &  1353 &  75 & 29&662&833&868  \\
192& Unaccelerated  &  1423 &2204 &46&267 &2483&2727 \\
\bottomrule
\end{tabular}
\end{table}

It is observed, from Table~\ref{tab:direct-vs-nda}, that  the memory 
usage in the unaccelerated transport solver is slightly 
higher than that in the diffusion 
accelerated transport method since the diffusion system 
takes a negligible amount of memory.  Take $N_\text{d}=32$
as an example,  where the memory consumed by  the diffusion 
accelerated transport method is less than that used in the unaccelerated
transport solver by $44$ M.  The memory  difference between 
the unaccelerated transport solver and the diffusion accelerated  transport  
algorithm  becomes larger  as more angular directions are used. 
The difference is $63$ M at $64$ angular directions, and it grows 
to $117$ M and $70$ M when $128$ and $192$ angular directions are 
employed.  The compute time on the matrix-free matrix-vector 
operations, in the diffusion accelerated transport solver, does not change 
much as more angular directions are added because the matrix-free matrix-vector 
operations occur in the low order diffusion system and the diffusion system stays
the same regardless of the number of angular directions. On the other hand,
the compute time of the matrix-free matrix-vector operations for the unaccelerated
transport solver is doubled as we double the number of angular directions since 
the transport problem becomes twice larger when the number of angular directions 
is doubled.  For instance,  the compute time spent on the matrix-free matrix-vector 
operations  is $246$ s when using $32$ angular directions, and it grows to 
$567$ s, $1241$ s  and $2204$ s  for 64, 128 and 192 angular directions, respectively.
The preconditioner setup time for the diffusion accelerated transport solver 
is slightly  less than that in the unaccelerated transport solver for 64, 128 and 192 
angular directions, and  it is almost  the same as that of the unaccelerated 
transport solver at 32 angular directions.   The preconditioner 
application time for the diffusion accelerated transport solver 
is three times as much as that in the unaccelerated transport solver regardless 
of the number of angular directions since much more GMRES iterations 
are required. However, the linear solver time of the diffusion accelerated  
transport algorithm is much less than that spent on the unaccelerated transport 
solver.  For instance,  at 32 angular directions, the unaccelerated transport 
solver takes $283$ s on GMRES, while the diffusion accelerated transport 
solver costs $179$ s on GMRES.  The unaccelerated transport solver uses twice compute time
as much as the diffusion accelerated transport method for 64 angular directions,
 and this ratio is increased to 3 for 128 and 192 angular directions. 
 For the overall simulation time, 
  the diffusion accelerated  transport solver is twice as fast as  the unaccelerated transport solver for 
  32 and 64 angular directions, and  it is three times  faster for 128 and 192 angular directions.
  We conclude that the nonlinear diffusion acceleration  technique together  with 
  SGMASM enhances  the transport criticality calculations 
  significantly.   
  
\begin{figure}
 \centering
 \includegraphics[width=0.49\linewidth]{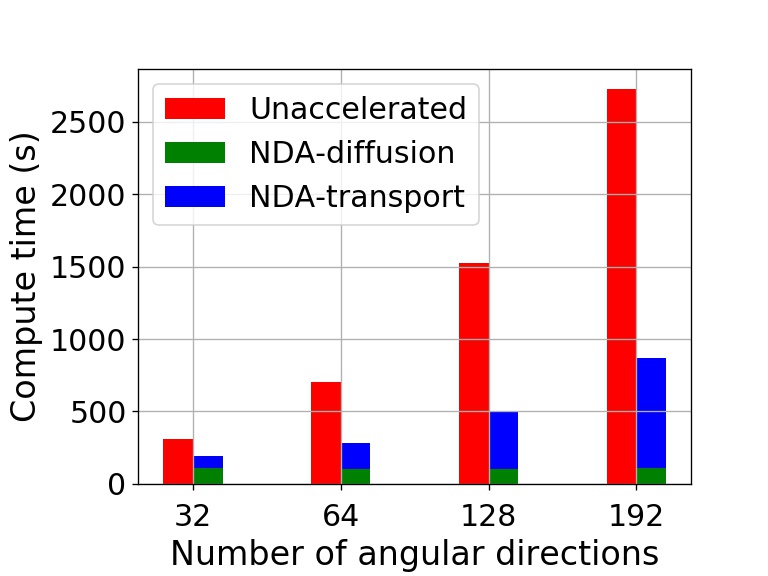} 
 \includegraphics[width=0.49\linewidth]{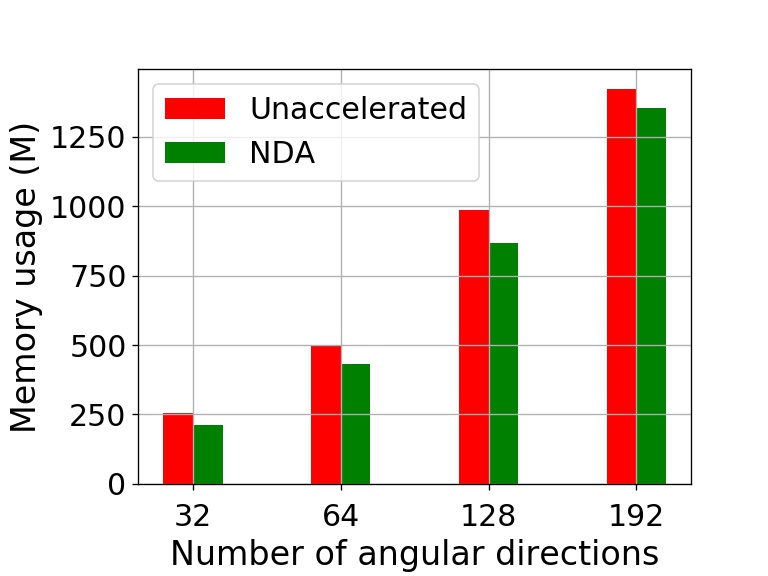} 
 \caption{Performance comparison in the compute time and the memory usage 
 with an unaccelerated transport solver using  different numbers 
 of angular directions.  Right: total compute time; left: memory usage.  \label{fig:direct-vs-nda}}
\end{figure}

  The performance comparison  in the compute time 
  and the memory usage is also drawn in Fig.~\ref{fig:direct-vs-nda}. 
Here ``NDA-diffusion" denotes the low order diffusion system, and ``NDA-transport" represents 
the transport system. It is obvious that the compute time used in the unaccelerated transport solver 
is much higher than that in the diffusion accelerated  transport  algorithm  especially when 
the number of angular directions  is large. For example, the diffusion accelerated transport solver 
is three times as fast as the unaccelerated transport solver when using 192 angular directions. 
The memory usages are similar to each other, and the 
memory usage in the unaccelerated  transport 
solver is slightly higher than that in the diffusion accelerated transport algorithm.

The compute times on the  individual components of the low-order diffusion system 
are also reported in Table~\ref{tab:low-order-problem-size} for 
different numbers of angular directions using 1,280 
processor cores. 
\begin{table}
\tiny
\centering
\caption{Compute times on the individual  components of the low-order diffusion 
system for different numbers of angular directions
 using 1,280  processor cores.  \label{tab:low-order-problem-size}}
\begin{tabular}{c c c c c c  c c c c}
\toprule
$N_\text{d}$  &Time$_{\text{ksp}}$ (s) & Time$_{\text{PCA}}$ (s) 
& Time$_{\text{PCS}}$ (s) &Time$_{\text{MF}}$ (s)
&Time$_{\text{Func}}$ (s)&Time$_\text{Jac}$ (s)& Time$_\text{LS}$ (s) & Time$_{\text{low}}$ (s) \\
\midrule
32 &95 &  15 &  6 & 77&4&4& 6 &106   \\
\midrule
64 &89 &  15 &  6 &71&4&  4&  6& 101  \\
\midrule
128 &95  &  15 &  6 & 72&4&4 &6&102 \\
\midrule
192& 97 &  20 &  6 & 75&4&5&6&109  \\
\bottomrule
\end{tabular}
\end{table}
As expected,  the compute times on the individual  components stay close to constants   
when we increase the number of angular directions. 

The trend does not hold  for the transport system as shown 
in Table~\ref{tab:problem-size-high-order}.
\begin{table}
\scriptsize
\centering
\caption{Compute times on the preconditioner setup and application of the transport system for different numbers of 
angular directions using 1,280 processor cores.  \label{tab:problem-size-high-order}}
\begin{tabular}{c c c c c c  c c c c}
\toprule
$N_\text{d}$  &Its$_\text{sweep}$ & Time$_{\text{PCS}}$ (s) & Time$_{\text{PCA}}$ (s) & Time$_\text{high}$ (s)  \\
\midrule
32 &35&4 &  73 & 88   \\
\midrule
64 &36&8 &  150 &   180  \\
\midrule
128 &39&14  &  338 & 400 \\
\midrule
192&39& 23&642  & 759  \\
\bottomrule
\end{tabular}
\end{table}
Here the averaged  number of GMRES iterations per Picard iteration 
stay close to a constant as more angular directions are used, which indicates that 
SGMASM is scalable in problem size.  
The compute times on the individual  components such as 
the preconditioner setup and the preconditioner application are doubled as expected  
when the number of angular directions are doubled. We conclude that the
diffusion accelerated transport solver  equipped with 
SGMASM is linearly  scalable in problem size.

\subsection{Influences of coarsening threshold}
In the matrix coarsening algorithm, there is a critical parameter  that   impacts not only the 
operator complexity  but also the convergence rate.  The parameter denoted as ``threshold"
determines   which elements are important to be kept in the strength matrix and which elements 
can be ignored.  More precisely, a connection from $i$ to $j$ is included  in the strength graph 
if and only if 
$$
-\vM(i,j) > \theta \max_{k \neq i} (-\vM(i,k)).
$$
Here $\vM(i,j)$ is the $j$th column of the $i$th row. 
If $\theta$ is too small, the operator complexity will be too high since all connections are considered. If 
$\theta$ is too large,  the algorithm will not  converge  since there are no enough coarse points 
to resolve low frequency   modes.  An optimal choice of $\theta$  is often problem dependent.  
In \cite{yang2002boomeramg},
$\theta=0.25$ is recommended for 2D elliptic problems and $\theta=0.5$  for 3D problems. 
The coarsening threshold is  introduced to 
construct  strength  matrices where only important coefficients  should be kept. 

\begin{figure}
 \centering
 \includegraphics[width=0.49\linewidth]{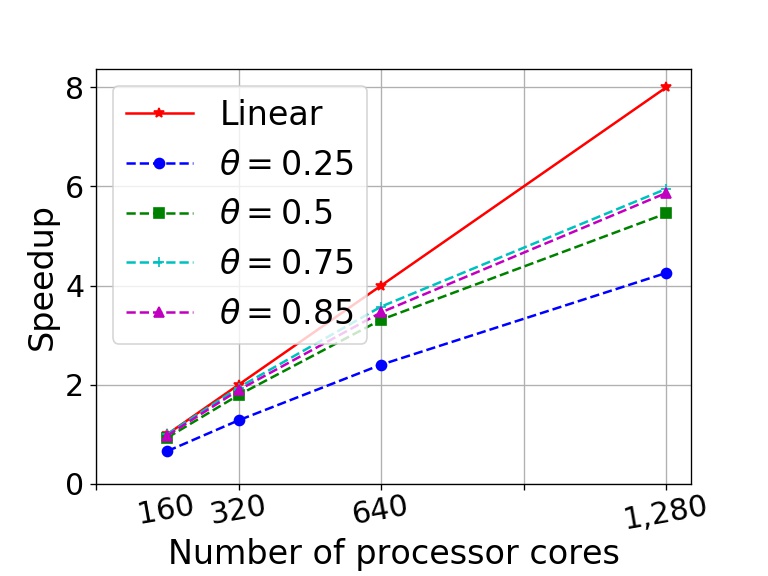} 
 \includegraphics[width=0.49\linewidth]{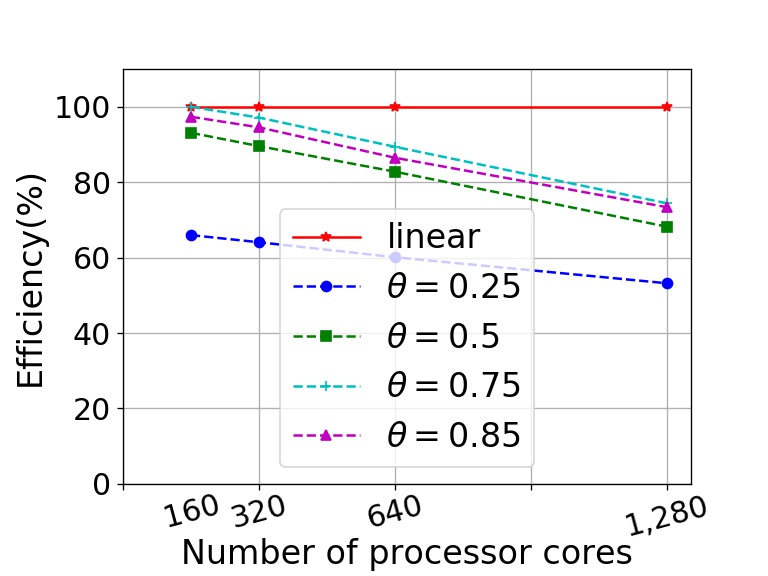} 
 \caption{Speedup and parallel efficiency for different coarsening thresholds  on
 160, 320, 640 and 1,280 processor cores.  \label{fig:threshold}}
\end{figure}

In this test, we study the algorithm performance with respect  
to $\theta$ on different numbers of processor cores. 
The same configuration as before is employed, and the numerical results are 
reported in Table~\ref{tab:threshold}. 
\begin{table}
\scriptsize
\centering
\caption{Impacts of coarsening threshold on the algorithm  performance using 
160, 320, 640 and 1,280 processor cores.  \label{tab:threshold}}
\begin{tabular}{c c c c c c  c c c c}
\toprule
$np$ & $\theta$ &Mem (M) & Its$_{\text{newton}}$ & Its$_{\text{linear}}$ 
&Its$_{\text{sweep}}$ &Time$_{\text{low}}$ (s)&Time$_\text{high}$ (s)& Time$_T$ (s)& EFF  \\
\midrule
160 &0.25 &  1587 &  2 & 32&30&663& 627 &1290 &66\%  \\
160 & 0.5  &  1526 &2 &25&19 & 445& 469& 914&93\% \\
160 & 0.75  &  1408 &2 &21&19& 452&  399& 851&100\% \\
160 & 0.85  &  1385 &2 &22&19& 462&  412& 874&97\% \\
\midrule
320 &0.25 &  810 &  2 & 33&30&350&  314&  664& 64\%  \\
320& 0.5 &  796 &2 &25&19 &234&  241& 475& 90\% \\
320 & 0.75  &  728 &2 &21&19& 236& 202& 438& 97\% \\
320 & 0.85  &  742 &2 &23&19& 241& 209& 450& 95\% \\
\midrule
640 &0.25  &  427 &  2 & 34&30&190 &164&354&  60\%  \\
640 & 0.5  &  408 &2 &26&19 & 129&128& 257& 83\% \\
640 & 0.75  &  380 &2 &22&19& 131&107& 238 &89\% \\
640 & 0.85  &  400 &2 &23&20& 134&112& 246 &86\% \\
\midrule
1,280& 0 .25&  210 &  2 & 35&31&111&89&200& 53\%  \\
1,280 & 0.5  &  216 &2 &27&19 &83&73&156& 68\% \\
1,280 & 0.75  &  211 &2 &23&19&83& 60&143&74\% \\
1,280 & 0.85  &  210 &2 &24&20&84& 61&145&73\% \\
\bottomrule
\end{tabular}
\end{table}
It is easily observed, from Table~\ref{tab:threshold}, 
that the number of Newton iterations stays 
as a constant for all   values of $\theta$ for all core counts. 
The memory usage often decreases as a larger $\theta$  is used.
For example, at 160 processor cores, the memory usage reduces from $1587$ M
to $1526$ M by  $61$ M when $\theta$ is increased from 0.25 to 0.5.
It continues being reduced to $1408$ M and $1385$ M,  when $\theta=0.75$ and 
$\theta=0.85$, respectively.   The same behaviors are observed  for high core counts
 as well; e.g., 
at 640 processor cores, the memory consumption  becomes  small  
to $408$ M and $380$ M from $427$ M as $\theta$ grows from 0.25 to 0.5 and 0.75.
It is because the operator complexity is  reduced  when a larger threshold is
employed.  There are some exceptions  where the memory usage is increased 
as a larger $\theta$ is used because  the averaged  number of GMRES iterations in 
the diffusion  system or the transport system or both becomes larger.  
The averaged  number of GMRES iterations in the diffusion 
system becomes smaller at the beginning when we increase $\theta$ from 
0.25 to 0.5, and then it does not change much any more when we
keep increasing $\theta$. For the 
160-core case, the averaged  number of GMRES iterations for the diffusion system 
at $\theta=0.25$  is 32, and it is reduced to 25 by 7  at $\theta=0.5$,  and to 
21 at $\theta=0.75$, but it is increased by one iteration at $\theta=0.85$, 
which indicates that $\theta=0.85$ is too big  for this particular problem.
Compared with the diffusion system, the threshold  has a bigger impact on 
the the transport system.  The  averaged number of GMRES iterations
per Picard iteration  is almost halved 
when $\theta=0.5$ is used instead of $\theta=0.25$; e.g., the number of GMRES
iterations is decreased  from 30 to 19 at 160 processor cores as 
$\theta$ is increased from 0.25 to 0.5. It is kept close to constants 
for $\theta=0.5, 0.75  \text{~and~} 0.85$.  The compute time in both 
the diffusion and the transport systems decreases significantly  at
the beginning, and then does not change much. This trend is 
consistent with the number of GMRES iterations. Take the 160-core case 
as an example, where the compute time in the diffusion system 
is reduced from 663 s to 445 s by 218 s when we increase the threshold from 0.25
to 0.5, and then it does not change much for $\theta=0.75$ and $0.85$. 
The transport  system has a similar trend, that is, 
it is decreased from 627 s to 469 s by 160 s  as $\theta$ is increased from 0.25 to 0.5 and
does not change much for $\theta=0.75$ and $\theta=0.85$.
The total simulation time has exactly the same trend as the compute times on the diffusion system 
and the transport system because it is a simple summation of the compute times of both. 
The parallel efficiency is computed using the smallest compute time obtained at 160 cores as a 
base so that  a higher efficiency represents a more efficient simulation. $\theta=0.25$
has a relatively low parallel efficiency, and all other choices have  good parallel 
efficiencies.  Even with up to 1,280 processor cores, the proposed algorithm 
is able to maintain a parallel efficiency above or around 70\%.  The proposed algorithm 
equipped with $\theta=0.75$ has the best performance, and it has a parallel efficiency 
 as high as $74\%$
on 1,280 processor cores.  The speedup and parallel efficiency
are also drawn in Fig.~\ref{fig:threshold}, where we observe  that $\theta=0.75$ is 
always better than the other choices.

To further  understand  the influences of coarsening thresholds, 
the compute times on the individual  components of 
the diffusion  system are summarized in Table~\ref{tab:threshold-low-order}. 
\begin{table}
\scriptsize
\centering
\caption{Compute times on the individual  components of the diffusion system 
with respect to different coarsening thresholds using 160, 320, 640, and 
1,280 processor cores.    \label{tab:threshold-low-order}}
\begin{tabular}{c c c c c c  c c c c}
\toprule
$np$ & $\theta$ &Time$_{\text{ksp}}$ 
& Time$_{\text{PCS}}$ (s) & Time$_{\text{PCA}}$ (s) &Time$_{\text{MF}}$ (s)
&Time$_{\text{Func}}$ (s) &Time$_\text{Jac}$ (s) & Time$_\text{LS}$ (s)  \\
\midrule
160 &0.25 &  585 &  15 & 39&548&28& 31 &37   \\
160 & 0.5  &  368 &11 &28&347 & 28& 31& 37 \\
160 & 0.75  &  374 &11 &30&371& 28&  31& 37 \\
160 & 0.85  &  384 &10 &30&362& 28&  31& 37 \\
\midrule
320 &0.25 &  310 &  9 & 21&289&15&  16&  19  \\
320& 0.5 &  194 &7 &16&180 &15&  16& 19 \\
320 & 0.75  &  196 &7 &17&181& 15& 16& 19 \\
320 & 0.85  &  201 &7 &17&186& 14& 16& 19 \\
\midrule
640 &0.25  &  170 &  7 & 15&152&7 &8&10 \\
640 & 0.5  &  108 &7 &13&93 & 8&8& 10 \\
640 & 0.75  &  110 &6 &13&95& 8&8& 10 \\
640 & 0.85  &  113 &5 &14&98& 8&8& 10 \\
\midrule
1,280& 0.25 &  100 &  6 & 10&86&4&4&5  \\
1,280 & 0.5  &  71 &8 &14&53 &4&4&6 \\
1,280 & 0.75  &  71 &6 &14&54&4& 4&6 \\
1,280 & 0.85  &  72 &7 &13&55&4& 5&6 \\
\bottomrule
\end{tabular}
\end{table}
It is found that the compute time on the linear solver decreases  significantly as $\theta$ is increased from 
0.25 to 0.5, and it stays close to constants for $\theta=0.75$ and $\theta=0.85$, which is consistent 
with the number of GMRES iterations. The preconditioner setup times are almost the same for 
all cases, and the preconditioner application time has a smilar trend as the linear solver time, that is,
it decreases  at the beginning and then does not changes  much.  As expected, the compute times 
on  the line search  operations, the function evaluations and the Jacobian evaluations 
stay the same for all values of $\theta$  because they are determined  
by the number  of Newton iterations that is a constant 
for all cases.  The compute time on the matrix-free matrix-vector operations  is completely  determined  by 
the number of GMRES iterations so that it has the same trend as the number of GMRES iterations. For instance,
at 640 processor cores,  the compute time of the matrix-free operations is reduced from 
152 s to 93 s by  $40\%$, and it does not change much for $\theta=0.75$ and $0.85$. $\theta=0.75$
is the best choice in the diffusion system in the sense that the compute time is the smallest.

 Similarly,  the compute times on the individual  components of the transport system are 
 also  drawn in Table~\ref{tab:threshold-high-order}. 
\begin{table}
\scriptsize
\centering
\caption{ Compute times on the individual  components of the transport system using  different 
coarsening thresholds on 160, 320, 640,  and 1,280 processor cores.   \label{tab:threshold-high-order}}
\begin{tabular}{c c c c c c  c c c c}
\toprule
$np$ & $\theta$ & Time$_{\text{PCS}}$ (s) & Time$_{\text{PCA}}$ (s) & Time$_\text{high}$ (s)  \\
\midrule
160 &0.25 &  18 &  549 & 627   \\
160 & 0.5  &  13 &417 &469 \\
160 & 0.75  &  12 &354 &399 \\
160 & 0.85  &  11 &365 &412 \\
\midrule
320 &0.25 &  10 &  273 & 314  \\
320& 0.5 &  8 &211 &241 \\
320 & 0.75  &  7 &177 &202 \\
320 & 0.85  &  7 &184 &209 \\
\midrule
640 &0.25  &  6 &  140 & 164 \\
640 & 0.5  &  5 &110 &128 \\
640 & 0.75  &  5 &93 &107 \\
640 & 0.75  &  4 &96 &112 \\
\midrule
1,280& 0.25 &  4 &  74 & 89  \\
1,280 & 0.5  &  5 &61 &73 \\
1,280 & 0.75  &  3 &51 &60 \\
1,280 & 0.85  &  3 &51 &61 \\
\bottomrule
\end{tabular}
\end{table}
Here the preconditioner setup time is almost the same for all coarsening 
thresholds  except at 160 processor cores,
where $\theta=0.25$ takes 50\% more compute time than that using other thresholds. 
The preconditioner application time  decreases a lot as $\theta$ is increased from 0.25 to 0.5, and it 
slightly decreases  again when using 
$\theta=0.75$, and then it becomes a little larger for $\theta=0.85$. 
For the 160-core case, the preconditioner application time  is
reduced from 549 s to 417 s by $24\%$ as $\theta$ grows from 0.25 to 0.5,  
and it continues  being  reduced to $354$ s by $15\%$ when  $\theta=0.75$,
and then it slightly increases  to $365$ s for $\theta=0.85$. This trend holds for all 
  core counts. The total compute time of the transport  system 
behaves  in the same way as the preconditioner  application time 
because the preconditioner application accounts  for $90\%$ of the total compute time.
Again,  $\theta=0.75$ is the best choice in terms of the compute time  for the transport system.

\subsection{Influences of the number of aggressive coarsening levels}
The complexity of coarse operators affects the algorithm performance  
in the memory usage and the compute time.  A high-complexity coarse 
operator uses more memory and costs more compute time per iteration.  
The complexity  of 
coarse operators can be reduced by introducing an aggressive coarsening scheme 
whose basic idea is  to keep as few coarse points as possible in the coarse levels 
\cite{kong2019highly, yang2002boomeramg}. 

In this test, 
we study the influences of the numbers of aggressive coarsening levels on the complexity 
of the coarse operators and on  the algorithm performance. The same configuration
as before is used. ``agg=0" denotes that no aggressive  coarsening is applied. 
The numerical results are summarized in Table~\ref{tab:agg}.
\begin{table}
\scriptsize
\centering
\caption{Influences  of  the number of aggressive coarsening levels 
on the algorithm performance 
using  160, 320, 640 and 1,280 processor cores.  \label{tab:agg}}
\begin{tabular}{c c c c c c  c c c c}
\toprule
$np$ & agg &Mem (M) &Comp & Its$_{\text{linear}}$ &Its$_{\text{sweep}}$
&Time$_{\text{low}}$ (s)&Time$_\text{high}$ (s)& Time$_T$ (s)& EFF  \\
\midrule
160 &0 &  2136 &  2.77 & 16&17&459& 508 &967 &88\%  \\
160 & 1  &  1627 &1.88 &20&18 & 449& 452& 901&94\% \\
160 & 2  &  1475 &1.64 &21&18& 446&  416& 862&99\% \\
160 & 4  &  1414 &1.55 &21&19& 454&  401& 862&99\% \\
160 & 8  &  1408 &1.55 &21&19& 451&  399& 850&100\% \\
\midrule
320 &0 &  1057 &  2.79 & 17&17&247&  263&  510& 83\%  \\
320& 1 &  864 &1.89 &20&18 &245&  230& 475& 89\% \\
320 & 2  &  755 &1.65 &21&18& 246& 213& 459& 93\% \\
320 & 4  &  730 &1.56 &22&19& 242& 203& 445& 96\% \\
320 & 8  &  728 &1.55 &22&19& 239& 200& 439& 97\% \\
\midrule
640 &0  &  542 &  2.83 & 17&17&141 &142&283&  75\%  \\
640 & 1  &  431 &1.92 &21&18 & 145&126& 271& 78\% \\
640 & 2  &  409 &1.66&22&18& 140&116& 256 &83\% \\
640 & 4  &  383 &1.57 &22&19& 135&109& 244 &87\% \\
640 & 8  &  380 &1.56 &22&19& 129&107& 236 &90\% \\
\midrule
1,280& 0&  303 &  2.86 & 18&17&99&81&180& 59\%  \\
1,280 & 1  &  239 &1.95 &22&18 &96&73&169& 63\% \\
1,280 & 2  &  230 &1.67 &23&18&89& 66&155&69\% \\
1,280 & 4  &  212 &1.57 &23&19&87& 62&149&71\% \\
1,280 & 8  &  211 &1.57 &23&19&80& 61&141&75\% \\
\bottomrule
\end{tabular}
\end{table}

It is found, from Table~\ref{tab:agg}, that  the memory usage becomes smaller 
as more aggressive  coarsening levels are employed. 
For the 160-core case, the memory usage  is reduced by $23\%$
from $2136$ M to $1627$ M when  one  aggressive coarsening level is introduced. 
It continues being  decreased to $1475$ M by $9\%$ and to $1414$
M by $4\%$ as 2 and 4 aggressive 
coarsening  levels are used,
respectively, and it stays almost the same for ``agg=8". 
The same pattern is observed for all other core counts, that is,
the memory usage is reduced significantly  for ``agg=1" and 
``agg=2", and stays almost the same for ``agg=4" and ``agg=8".
This pattern occurs because of the complexity of the coarse
operators.  The complexity drops much when ``agg=1" is employed
instead of ``agg=0".  For example, at 1,280 processor cores,
the complexity is reduced from 2.86 to 1.95 by $31\%$  when ``agg=1" is 
used, and continues becoming  a little smaller as ``agg=2", ``agg=4" and 
``agg=8" are adopted. 
We do not list the number of Newton iterations since it is kept as a constant, 2, 
for all cases. The number of GMRES iterations in the diffusion  system  becomes 
larger as more aggressive coarsening  levels are used,
but the performance of  the resulting algorithm 
does not deteriorate because the per-iteration cost decreases 
significantly.   The number of  GMRES iterations in 
the transport system stays  close to constants for different numbers 
of aggressive coarsening levels. 
The number of aggressive  coarsening levels 
does not influence  the compute time  much in the diffusion system 
when the number of processor cores is small, but it has 
a bigger impact   when we use more  cores. 
For example, the compute time of the diffusion system is reduced to 80 s from 99
s by $20\%$ when ``agg=8" is employed.
On the other hand, the number of aggressive coarsening levels 
has a consistent impact on the algorithm performance for the transport system. 
More precisely, as ``agg=8" is employed the compute time is reduced by $25\%$ -
$30\%$.  For example, for the 320-core case, 
the compute time in the transport  system reduces  from $263$ s to $230$ s,
$213$ s, $203$ s and $200$ s when the number of aggressive coarsening levels 
is increased from 0 to 1, 2, 4 and 8.  The compute time at `agg=8` is  $24\%$
less than that using ``agg=0".  The total compute  time correlates to  the
complexity  of the coarse operators, that is, it decreases at the beginning   and then does
not change much.  The parallel efficiency is computed using the smallest compute time
obtained with 160 processor cores as a base. The proposed
 algorithm with ``agg=8" has the best parallel efficiency, and a 
 parallel efficiency as high as $75\%$ is obtained with up to 1,280 processor cores.  
 
\begin{figure}
 \centering
 \includegraphics[width=0.49\linewidth]{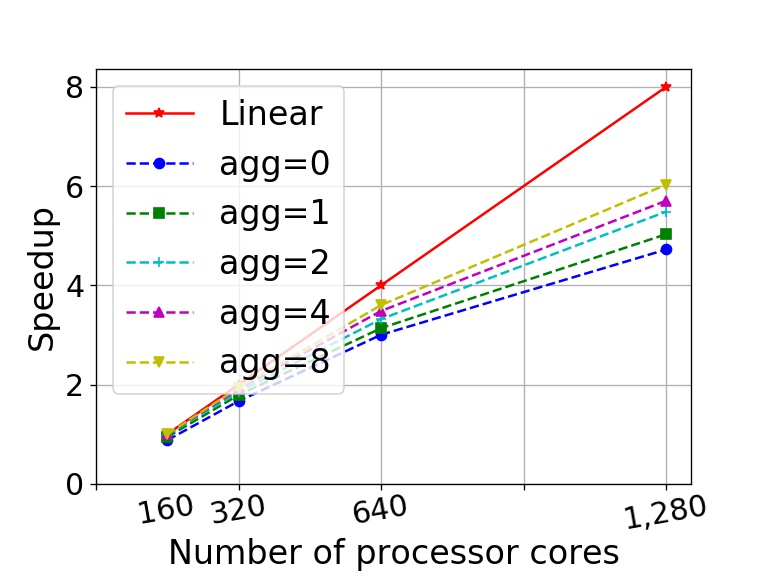} 
 \includegraphics[width=0.49\linewidth]{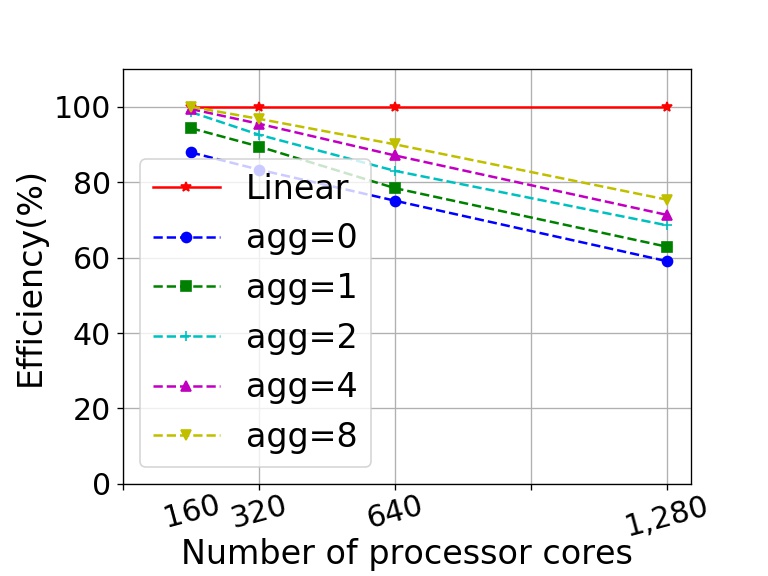} 
 \caption{Speedup and parallel efficiency for different numbers of aggressive coarsening levels   on
 160, 320, 640 and 1,280 processor cores. Left: speedup; right: parallel efficiency.  \label{fig:agg}}
\end{figure}

 The speedup 
 and the parallel efficiency are also drawn in Fig.~\ref{fig:agg}, where we 
found  that a larger  number of aggressive  coarsening levels often leads to a better parallel efficiency 
for all core counts. 

We next study the impacts of the number of aggressive coarse levels on the individual   
components of the diffusion system.  The compute times on the individual  
components of the diffusion system using  different numbers of aggressive coarse levels 
are reported in Table~\ref{tab:agg-low-order}. 
\begin{table}
\scriptsize
\centering
\caption{Compute times on the individual   components of the diffusion system 
with respect to different numbers of aggressive coarsening levels using 160, 320, 640
and 1280 processor cores.    \label{tab:agg-low-order}}
\begin{tabular}{c c c c c c  c c c c}
\toprule
$np$ & agg&Time$_{\text{ksp}}$ (s)& Time$_{\text{PCS}}$ (s)
& Time$_{\text{PCA}}$ (s) &Time$_{\text{MF}}$ (s)
&Time$_{\text{Func}}$ (s)&Time$_\text{Jac}$ (s)& Time$_\text{LS}$ (s)  \\
\midrule
160 &0 &  380 &  25 & 53&320&28& 31 &37   \\
160 & 1  &  371 &17 &40&332 & 28& 31& 37 \\
160 & 2  &  369 &14 &38&335& 28&  31& 37 \\
160 & 4  &  376 &12 &34&348& 28&  31& 37 \\
160 & 8  &  373 &11 &30&350& 28&  31& 37 \\
\midrule
320 &0 &  206 &  17 & 35&163&15&  16&  19  \\
320& 1 &  245 &11 &29&173 &14&  16& 19 \\
320 & 2  &  206 &10 &27&177& 15& 16& 19 \\
320 & 4  &  202 &8 &22&181& 15& 16& 19 \\
320 & 8  &  198 &8 &18&182& 14& 17& 19 \\
\midrule
640 &0  &  120 &  12 & 27&85&8 &8&10 \\
640 & 1  &  123 &9 &25&93 & 8&9& 10 \\
640 & 2  &  119 &8 &23&92& 8&8& 10 \\
640 & 4  &  114 &6 &18&94& 7&8& 10 \\
640 & 8  &  108 &6 &12&94& 8&8& 10 \\
\midrule
1,280& 0 &  88 &  11 & 29&49&4&5&6  \\
1,280 & 1  &  85 &9 &27&50 &4&4&6 \\
1,280 & 2  &  78 &7 &22&51&4& 4&6 \\
1,280 & 4  &  75 &7 &17&54&4& 5&6 \\
1,280 & 8  &  69 &6 &12&54&4& 4&5 \\
\bottomrule
\end{tabular}
\end{table}
As more processor cores are used, 
the impacts of the number of aggressive coarsening  levels on the linear solver time 
become larger.  When using 160 processor cores, the compute time on 
the linear solver in the diffusion system does not change much with an increase in 
the number of aggressive coarsening levels, while the linear solver time obtained using ``agg=8"
is $22\%$ less than that with ``agg=0" when the number of processor cores is
1,280.  The preconditioner setup time is reduced significantly when the number of aggressive coarsening levels
is increased from 0 to 1, and then continues gradually decreasing   for ``agg=2", `agg=4" and `agg=8".
Take the 1,280-core case as an example, where the preconditoner setup time 
is reduced by $50\%$ from $11$ s to $6$ s when ``agg=8" is 
employed instead of ``agg=0". This pattern is shown for the preconditioner application time as well.
It decreases  much when one  aggressive coarsening level is adopted, and slightly 
reduces as more aggressive levels are used. For example, for the 320-core case,
the preconditioner application time reduces  to 29 s from 35 s when the number of 
aggressive coarsening levels is increased from 0 to 1, and continues being 
decreased to $27$ s , $22$ s and $18$ s when 
the number of aggressive coarsening  levels is increased to 2, 4  and 8.
The compute time on the matrix-free operations is increased as more aggressive  coarsening levels 
are used because the resulting algorithm has more GMRES iterations. 
However, the time increase on the matrix-free operations  is compensated by  the time reduction 
on the precontioner setup and application so that the total compute time 
is properly  decreased  as more aggressive coarsening levels are employed.
The compute times on the Jacobian evaluations, the function evaluations and the line search 
are almost exactly the same for all cases since they are determined by
the number of Newton iterations that is the same for all tests.

Similarly,  we report the compute times on the individual 
  components  of  the transport system in Table~\ref{tab:agg-high-order}.
\begin{table}
\scriptsize
\centering
\caption{Compute times on the individual   components of the transport system using different
numbers of aggressive  coarsening levels on 160, 320, 640 and 1,280 processor cores.  \label{tab:agg-high-order}}
\begin{tabular}{c c c c c c  c c c c}
\toprule
$np$ & agg & Time$_{\text{PCS}}$ (s) & Time$_{\text{PCA}}$ (s) & Time$_\text{high}$ (s)  \\
\midrule
160 &0 &  28 &  457 & 508   \\
160 & 1  &  18 &403 &452 \\
160 & 2  &  14 &369 &416 \\
160 & 4  &  12 &356 &401 \\
160 & 8  &  12 &353 &399 \\
\midrule
320 &0 &  17 &  233 & 263  \\
320& 1 &  11 &202 &230 \\
320 & 2  &  9 &187 &213 \\
320 & 4  &  7 &179 &203 \\
320 & 8  &  7 &178 &200 \\
\midrule
640 &0  &  11 &  123 & 142 \\
640 & 1  &  7 &110 &126 \\
640 & 2  &  6 &100 &116 \\
640 & 4  &  5 &94 &109 \\
640 & 8 &  5 &93 &107 \\
\midrule
1,280& 0 &  7 &  70 & 81  \\
1,280 & 1  &  6 &61 &73 \\
1,280 & 2  &  4 &56 &66 \\
1,280 & 4  &  3 &53 &62 \\
1,280 & 8  &  4 &50 &61 \\
\bottomrule
\end{tabular}
\end{table}
The same pattern observed earlier  presents  in  the transport system. 
More precisely, the preconditioner setup time for the transport system 
decreases much at the beginning, and then gradually  deceases as 
more aggressive coarsening levels are used.  For example,
at 640 processor cores, the preconditioner setup takes 11 s 
when ``agg=0", and the compute time is halved  to 5 s when the number of 
aggressive coarsening levels is 8.  Similarly,  the precontioner application time
and the total compute time decrease as we employ more   aggressive coarsening levels.  
We conclude that in this test ``agg=8"
is the best in terms of the memory usage and the compute time for the overall simulation.
A large number of aggressive  coarsening  levels are required to maintain a good 
scalability.

\subsection{Comparison with a traditional multilevel method}
In this test, we compare the performance of the monolithic  
multilevel method equipped with the subspace-
based coarsening algorithm   with that of the 
unmodified  traditional  multilevel method using the full space 
based  coarsening approach.
The subspace-based coarsening algorithm mainly aims 
at reducing the preconditioner setup cost, and meanwhile 
is able to improve the preconditioner application performance.
The mesh used in this test is finer than that in the previous tests  since more 
processor cores will be used. The mesh has 6,464,825
nodes and 12,543,552 elements. The diffusion system has 45,253,775 unknowns,
and the transport system has  2,896,241,600 unknowns with 64 angular 
directions. The numerical results are summarized in Table~\ref{tab:fullspace}.
\begin{table}
\scriptsize
\centering
\caption{Performance comparison btweenn SGMASM  and  MASM using up to 5,120 processor cores.  \label{tab:fullspace}}
\begin{tabular}{c c c c c c  c c c c}
\toprule
$np$ & Algorithm & Time$_{\text{PCS}}$ (s)
& Time$_{\text{PCA}}$ (s)  &Time$_{\text{MF}}$ (s)
&Time$_{\text{Func}}$ (s)& Time$_\text{T}$ (s) & EFF \\
\midrule
1,280 &SGMASM &  47 &  1171 & 391 & 29 & 1802 & 100\%   \\
1,280 & MASM  &  -- &-- &-- & -- & -- & -- \\
\midrule
2,560 & SGMASM&  41 &  751 & 207 & 15 & 1105 & 82\%  \\
2,560& MASM &  401 &1020 &250 & 19  & 1700 & 53\%\\
\midrule
5,120 &SGMASM  &  28 &  476 & 160 & 11 & 734 & 61\% \\
5,120 & MASM  &  251 &920  & 136 & 10 & 1362 & 33\% \\
\bottomrule
\end{tabular}
\end{table}

It is found that the precontioner setup time is significantly reduced by using 
SGMASM.  The data misses for MASM at 1,280 processor cores 
because MASM requires the amount of  memory beyond the machine memory limit.
For 2,560 and 5,120 processor cores, the preconditioner setup of SGMASM
is ten times as fast as that of MASM. 
The preconditioner setup of MASM takes $401$ s at 2,560 processor cores and $251$ s at 
5,120 processor cores, while that of SGMASM 
costs only $41$ s and $28$ s, respectively. 
The preconditioner  application is improved by a factor of 2, compared with MASM. 
The preconditioner  application  time obtained using 
MASM is twice as high as that used by SGMASM, that is,
SGMASM is once  faster than MASM in the preconditioner  application for 2,560 and 
5,120 processor cores.
The compute times on the matrix-free operations  are similar  to each other for all core counts. 
Similarly, the compute times on the function evaluations  are close to each other.
For the overall simulation, SGMASM is able to run 
once faster than MASM.  SGMASM is capable of maintaining  a good 
parallel  efficiency  with up to 5,120 processor cores, while MASM is inefficient.

\section{Concluding remarks}
A nonlinear diffusion acceleration method  has been studied to improve 
the transport criticality  calculations, where the scattering and 
fission terms are evaluated using the computed scalar fluxes and 
eigenvalue from the nonlinear diffusions equations. 
To compute the  eigenvalue  of the low order diffusion system,
an inexact Jacobian-free Newton with a few  inverse power  iterations
as an initial  guess was   employed, and during each Newton iteration
a parallel monolithic multilevel preconditioner 
  together with GMRES was  adopted for calculating  
the Jacobian system.  The monolithic multilevel   method was
 also used for the solution of 
the linear system of transport equations.
To reduce the cost on the coarse space construction of the multilevel  method,
we studied  a subspace-based coarsening algorithm  that has been shown  to be 
more efficient than a traditional  full-space 
coarsening  approach on  thousands   of processor cores for an unstructured 
mesh problem with billions of unknowns. 
We have numerically verified  that 
the overall algorithm equipped with several  important ingredients; e.g., subspace-based coarsening,
monolithic  coupling, strength matrix  thresholding and aggressive coarsening, 
is scalable with up to  thousands of processor cores. 

While this work focuses on a multilevel Schwarz preconditioner
 for accelerating neutron transport calculations,
other   techniques such as Method of Characteristics (MOC) 
\cite{shaner2016verification, yamamoto2017genesis} may be
explored in the future.  Only strong scaling results were 
reported in this work, and weak scaling will be explored in our future work.

\section*{Acknowledgments}
This manuscript has been authored by Battelle Energy Alliance,
LLC under Contract No. DE-AC07-05ID14517 with the U.S. Department of Energy.
The United States Government retains and the publisher,
by accepting the article for publication, acknowledges
that the United States Government retains a nonexclusive,
paid-up, irrevocable, world-wide license to publish or
reproduce the published form of this manuscript,
or allow others to do so, for United States Government purposes.

This research made use of the resources of the High-Performance 
Computing Center at Idaho National Laboratory, which is supported by 
the Office of Nuclear Energy of the U.S. Department of Energy and the 
Nuclear Science User Facilities under Contract No. DE-AC07-05ID14517.  

\section*{References}

\bibliographystyle{elsarticle-num}

\bibliography{ndascaling}

\end{document}